# Piecewise Probability Distribution Theory

## by

## Lev Gelimson

**The "Collegium" All World Academy of Sciences**

**Munich (Germany)**

2nd Edition (2013)

1st Edition (2012)

## Abstract

A general piecewise (including pointwise) probability distribution with space-saving notation and its hierarchical particular cases are considered. The explicit closed-form normalization, expectation, and variance formulas along with the median and mode formulas and algorithms for a general one-dimensional piecewise linear probability distribution are obtained. They are also applied to a general polygonal, or one-dimensional piecewise linear continuous, probability distribution and, in particular, to a tetragonal probability distribution. The known formulas for the last distribution and a triangular probability distribution as a further particular case are used to test the obtained formulas and algorithms.

# Introduction

Both particular and some more general mostly continuous (continual without discontinuity points and jumps) piecewise linear probability distributions which can also be multidimensional are well known [Cramér]. For a triangular probability distribution, some basic formulas are also well known [Kotz & van Dorp, Wikipedia Triangular distribution]. [Kotz & van Dorp, van Dorp & Kotz] introduced trapezoidal distributions also generalized including partially nonlinear. [Kim] took into account truncated triangular and trapezoidal distributions also with discontinuity points and jumps. [Karlis & Xekalaki] considered such polygonal distributions only which are representable via positively weighted sums of independent triangular probability distributions on the unit segment [0, 1]. The mean of such a known polygonal distribution is placed in the middle third part of this unit segment only.

The present work is devoted to presenting a general piecewise (including pointwise) probability distribution with space-saving notation and its hierarchical particular cases. The same holds for analytically solving some fundamental problems for piecewise linear probability distributions which may be discontinuous and are namely directly introduced, which ensures most possible generality. They are very simple, natural, and typical and can provide adequately modeling via efficiently approximating practically arbitrary nonlinear probability distributions with any desired or/and required precision. General polygonal, or one-dimensional piecewise linear continuous, probability distributions are also very important extensions of tetragonal and triangular probability distributions. It is very natural to verify analytical methods of solving problems for general piecewise linear probability distributions via using some well-known basic formulas for a generalized trapezoidal probability distribution and for a triangular probability distribution. Geometrical approach can be also used to additionally verify analytical methods. If there are too many possible cases, which is typical for any piecewise problems, then apply algorithmic approach rather than explicit analytical closed-form solutions. The problems of the existence and uniqueness of the mean, median, and mode values for a general one-dimensional piecewise linear probability distribution are often nontrivial and can be of great importance for practice. It is very useful to provide clear mathematical (probabilistic and statistical) sense of methods and results. Setting and solving many typical urgent problems is the only criterion of creating, developing, and estimating any new useful theory. There are such problems not only in probability theory and mathematical statistics, but also in physics, engineering, chemistry, biology, medicine, geology, astronomy, meteorology, agriculture, politics, management, economics, finance, psychology, etc.

# 1. General Piecewise (Including Pointwise) Probability Distribution with Space-Saving Notation

## 1.1. General Piecewise Function

Consider a general piecewise function. Using space-saving notation [Gelimson 2012a], piecewise represent any function $g_{Ra}(x_D)$ on domain D by $x \in D$ with range Ra as a domain of dependent variable g by value $g(x) \in Ra$. Use subfunctions $g_{Ra(j)}(x_{D(j)})$ on non-intersecting subdomains $D_j$ by $x \in D_j$ with ranges $Ra_j$ as follows:

$$g_{Ra}(x_D) = \cup_{j \in J}\, g_{Ra(j)}(x_{D(j)}).$$

Here
symbol $\cup$ unifies subfunctions on subdomains similarly to set theory and can be also indexed with an index set and range,
J is any (possibly uncountable) index set and range,

$$D(j) = D_j$$

are non-intersecting subdomains of domain D so that

$$D = \cup_{j \in J}\, D_j$$

whereas range Ra is the union of all its subranges $Ra_j$ :

$$Ra = \cup_{j \in J}\, Ra_j\,.$$

Notata bene:
1. A partition, or non-intersecting distribution, of a domain between its subdomains is theoretically preferable.
2. However, it is possible that at a common point of continuity (e.g. at a boundary point) x of some subdomains

$$D_j = D(j) \ni x$$

with some subset

$$J_x \subseteq J$$

of indexes

$$j \in J_x\,,$$

all the partial values $g_j(x)$ coincide and hence build common value $g(x)$. Then it is admissible to explicitly include such point x into all these subdomains

$$D_j = D(j) \mid j \in J_x\,.$$

3. It is also possible that at a common point of continuity (e.g. at a boundary point) x of some subdomains

$$D_j = D(j) \ni x$$

with some subset

$$J_x \subseteq J$$

of indexes

$$j \in J_x\,,$$

some partial values $g_j(x)$ for

$$j \in J_{x|g=0}$$

vanish whereas all the remaining partial values $g_j(x)$ for

$$j \in J_{x|g \neq 0}$$

coincide and hence build common nonzero value $g(x)$. Then it is admissible to explicitly include such point x into all the subdomains

$$D_j = D(j) \mid j \in J_{x|g \neq 0}\,.$$

Support
$$S = \mathrm{supp}(g(x))$$
of a function $g(x)$ is the set of all x (from domain D) for which $g(x)$ is namely nonzero:
$$S = \mathrm{supp}(g(x)) = \{x \in D \mid g(x) \neq 0\},$$
$$S \subseteq D.$$
Support S may be conditionally extended via weakening its requirement. Namely, require that
$$g(x) = 0 \mid x \in D \setminus S$$
rather than
$$g(x) \neq 0 \mid x \in S$$
so that $g(x)$ may vanish at some points $x \in S$.

Compact support [S] of a function $g(x)$ is the smallest compact extension of support S, i.e. the intersection of all the compact extensions of support S.

## 1.2. General Pointwise Function

Consider a pointwise function as a particular case of a piecewise function. Regard all the separate distinct elements of a domain (as a set) as its subdomains (subsets). Identify [Gelimson 2003a, 2003b] one-point set {x} at least here with this element (point) x itself. Use this element as an index and a whole domain as an index set. Then a pointwise function is as follows:

$$g_{Ra}(x_D) = \cup_{x \in D} \, g_{\{g(x)\}}(x_{\{x\}}),$$

or, simplifying,

$$g_{Ra}(x_D) = \cup_{x \in D} \, g_{g(x)}(x_x),$$

or simply

$$g_{Ra}(x_D) = \cup_{x \in D} \, g(x).$$

Here obvious domains $\{g(x)\} = g(x)$ and $\{x\} = x$ of variables g and x, respectively, can be omitted, D is any (possibly uncountable) index set, all one-element sets {x} are subdomains of domain D so that

$$D = \cup_{x \in D} \, \{x\}$$

whereas range Ra is the union of all its subranges $Ra_j$ :

$$Ra = \cup_{x \in D} \, \{g(x)\}.$$

## 1.3. General Piecewise Probability Distribution

To consider namely a probability density function (distribution) f(x), take a non-negative real-valued function g(x), in our case a non-negative real-valued piecewise function
$$g_{Ra}(x_D) = \cup_{j \in J} g_{Ra(j)}(x_{D(j)}).$$
Here range Ra and all its subranges $Ra_j$ are subsets of the set
$$R_0^+ = [0, +\infty)$$
of all the non-negative real numbers:
$$g_{Ra \subseteq [0, +\infty)}(x_D) = \cup_{j \in J} g_{Ra(j) \subseteq [0, +\infty)}(x_{D(j)}).$$
Further in order to provide a non-negative real-valued function $f_{Ra}(x_D)$ with the role of a probability density function, the integral normalization condition
$$\int_D f(x) \, dD = 1$$
has to be satisfied. Each of these both conditions (of non-negativity and normalization) are necessary, and their pair is sufficient for a possibility of f(x) to be a probability density function. Then, beginning with a non-negative real-valued piecewise function
$$g_{Ra}(x_D) = \cup_{j \in J} g_{Ra(j)}(x_{D(j)})$$
with namely positive integral
$$\int_D g(x) \, dD > 0,$$
simply divide this function g(x) by this integral to obtain a probability density function:
$$f_{Ra}(x_D) = g_{Ra}(x_D) / \int_D g_{Ra}(x_D) \, dD = \cup_{j \in J} g_{Ra(j)}(x_{D(j)}) / \Sigma_{j \in J} \int_{D(j)} g_{Ra(j)}(x_{D(j)}) \, dD_j \, .$$
Support S of a non-negative real-valued function g(x) is the set of all x (from domain D) for which g(x) is namely strictly positive:
$$S = \{x \in D \mid g(x) > 0\},$$
$$S \subseteq D \, .$$
Compact support [S] of a non-negative real-valued function g(x) is the smallest compact extension of support S, i.e. the intersection of all the compact extensions of support S.
By integration, we may simply replace domain D and subdomains $D_j$:
1) either via compact support [S] and compact subsupports $[S_j]$, respectively,
2) or via support S and subsupports $S_j$, respectively.
Then we have, e.g.,
$$\int_S f(x) \, dS = 1,$$
$$\int_{[S]} f(x) \, d[S] = 1,$$
$$\int_S g(x) \, dS > 0,$$
$$\int_{[S]} g(x) \, d[S] > 0,$$
$$f_{Ra}(x_D) = g_{Ra}(x_D) / \int_S g_{Ra}(x_D) \, dS = \cup_{j \in J} g_{Ra(j)}(x_{D(j)}) / \Sigma_{j \in J} \int_{S(j)} g_{Ra(j)}(x_{D(j)}) \, dS_j \, ,$$
$$f_{Ra}(x_D) = g_{Ra}(x_D) / \int_{[S]} g_{Ra}(x_D) \, d[S] = \cup_{j \in J} g_{Ra(j)}(x_{D(j)}) / \Sigma_{j \in J} \int_{[S(j)]} g_{Ra(j)}(x_{D(j)}) \, d[S_j].$$

## 1.4. General Pointwise Probability Distribution

To consider namely a probability density function (distribution) f(x), take a non-negative real-valued function, in our case a non-negative real-valued pointwise function

$$g_{Ra}(x_D) = \cup_{x \in D} g_{\{g(x)\}}(x_{\{x\}}),$$

or, simplifying,

$$g_{Ra}(x_D) = \cup_{x \in D} g_{g(x)}(x_x),$$

or simply

$$g_{Ra}(x_D) = \cup_{x \in D} g(x).$$

Consider a pointwise function as a particular case of a piecewise function. Regard all the separate distinct elements of a domain (as a set) as its subdomains (subsets). Identify [Gelimson 2003a, 2003b] one-point set {x} at least here with this element (point) x itself. Use this element as an index and a whole domain as an index set.

Here obvious domains {g(x)} = g(x) and {x} = x of variables g and x, respectively, can be omitted, D is any (possibly uncountable) index set, all one-element sets {x} are subdomains of domain D so that

$$D = \cup_{x \in D} \{x\}$$

whereas range Ra is the union of all its subranges {g(x)}:

$$Ra = \cup_{x \in D} \{g(x)\}.$$

Range Ra and all its subranges {g(x)} are subsets of the set

$$R_0^+ = [0, +\infty)$$

of all the non-negative real numbers:

$$g_{Ra \subseteq [0, +\infty)}(x_D) = \cup_{j \in J} g_{\{g(x)\} \subseteq [0, +\infty)}(x_{\{x\}}).$$

Further in order to provide a non-negative real-valued function $f_{Ra}(x_D)$ with the role of a probability density function, the integral normalization condition

$$\int_D f(x) \, dD = 1$$

has to be satisfied. Each of these both conditions (of non-negativity and normalization) are necessary, and their pair is sufficient for a possibility of f(x) to be a probability density function. Then, beginning with a non-negative real-valued piecewise function

$$g_{Ra}(x_D) = \cup_{x \in D} g_{\{g(x)\}}(x_{\{x\}})$$

with namely positive integral

$$\int_D g(x) \, dD > 0,$$

simply divide this function g(x) by this integral to obtain a probability density function:

$$f_{Ra}(x_D) = g_{Ra}(x_D)/\int_D g_{Ra}(x_D) \, dD = \cup_{x \in D} g_{\{g(x)\}}(x_{\{x\}})/\int_D g_{Ra}(x_D) \, dD = \cup_{x \in D} g(x)/\int_D g(x) \, dD \, .$$

Nota bene: Here unions may be also uncountable.
By integration, we may simply replace domain D and subdomains $D_j$:
1) either via compact support [S] and compact subsupports $[S_j]$, respectively,
2) or via support S and subsupports $S_j$, respectively.
Then we have, e.g.,

$$\int_S f(x) \, dS = 1,$$
$$\int_{[S]} f(x) \, d[S] = 1,$$
$$\int_S g(x) \, dS > 0,$$
$$\int_{[S]} g(x) \, d[S] > 0,$$
$$f_{Ra}(x_D) = g_{Ra}(x_D)/\int_S g_{Ra}(x_D) \, dS = \cup_{x \in D} g_{\{g(x)\}}(x_{\{x\}})/\int_S g_{Ra}(x_D) \, dS = \cup_{x \in D} g(x)/\int_S g(x) \, dS \, ,$$
$$f_{Ra}(x_D) = g_{Ra}(x_D)/\int_{[S]} g_{Ra}(x_D) \, d[S] = \cup_{x \in D} g_{\{g(x)\}}(x_{\{x\}})/\int_{[S]} g_{Ra}(x_D) \, d[S] = \cup_{x \in D} g(x)/\int_{[S]} g(x) \, d[S].$$

## 1.5. Particular Case Hierarchy

Domain D , support S , compact support [S], their partitions into subdomains $D_j$ , subsupports $S_j$ , compact subsupports $[S_j]$ , respectively, as well as a non-negative real-valued piecewise function g(x) with namely positive integral on D , may be arbitrary. Therefore, we have here a multidimensional hierarchy of particular cases.

1. Domain D may be, e.g., discrete, continual, or mixed (with both discrete and continual parts). These simplest possibilities are especially typical in practice. Domain D may be, in particular, a subset of a countably dimensional space or of finitely dimensional Euclidean space $R^n$ (n ∈ N = {1, 2, ...}), e.g. of one-dimensional Euclidean space R = (-∞ , ∞). Naturally, continual domain D may be, in particular, one of these spaces as a whole.

2. Support S may be any subset of domain D .

3. Compact support [S] is the smallest compact extension of support S , i.e. the intersection of all the compact extensions of support S . The Bolzano-Weierstrass theorem [Encyclopaedia of Mathematics] proved that in any Euclidean space, a set is compact if and only if it is both closed and bounded. Then a compact support coincides with the corresponding closed support, or the support closure.

4. Partition of domain D into subdomains $D_j$ may be arbitrary: from using domain D itself with no partition to separating every point x ∈ D . In any Euclidean space, partitioning domain D by every coordinate into a finite set of non-intersecting intervals at least partially containing their endpoints is typical if possible. These simplest of the additive Borel sets provide the identity of all the common measures [Cramér, Encyclopaedia of Mathematics] and are preferable in probability theory, too.

5. Partition of support S into subsupports $S_j$ may be arbitrary: from using support S itself with no partition to separating every point x ∈ S . In any Euclidean space, partitioning namely bounded support S by every coordinate into a finite set of non-intersecting intervals at least partially containing their endpoints is typical if possible. These simplest of the additive Borel sets provide the identity of all the common measures [Cramér, Encyclopaedia of Mathematics] and are preferable in probability theory, too.

6. Partition of compact support [S] into subsupports $[S_j]$ may be arbitrary: from using support [S] itself with no partition to separating every point x ∈ [S] . In any Euclidean space, partitioning namely bounded compact support [S] by every coordinate into a finite set of intervals at least partially containing their endpoints and intersecting at them is typical if possible. These simplest of the additive Borel sets provide the identity of all the common measures [Cramér, Encyclopaedia of Mathematics] and are preferable in probability theory, too.

7. A non-negative real-valued piecewise function g(x) with namely positive integral on D may be arbitrary. To provide integrability of probability density function f(x) also multiplied by desired and/or required powers of variable x to obtain explicit (closed-form) integral (cumulative) probability distribution function F(x) along with moments [Cramér, Encyclopaedia of Mathematics], use namely the simplest and most suitable classes of functions to piecewise build a desired and/or required non-negative real-valued function g(x). Among them are, e.g., some power functions including polynomials, rational, exponential, trigonometric, and hyperbolic functions, as well as their linear and nonlinear combinations. Such function variety and partition variety provide very many possibilities of solving typical classes of urgent problems also in probability theory and mathematical statistics.

Therefore, namely the simplest piecewise linear functions whose domain D is a one-dimensional space R = (-∞ , ∞) and whose support S is bounded and representable via a finite set of non-intersecting intervals at least partially containing their endpoints are further used in the present work. They provide adequately fitting practically any desired and/or required function.

## 1.6. Nonstrictly Monotonic Function Continualization and Inversion

Strictly monotonic function inversion is well-known [Encyclopaedia of Mathematics] and straightforward because every image (output) always has namely the only preimage (input). But the common notation of the inverse is not suitable in general. For example, the inverse to function
$$y = h(x)$$
is usually denoted by
$$x = h^{-1}(y)$$
with always necessary explanation that $h^{-1}(y)$ means here NOT
$$1/h(y),$$
which would be natural, but the function which is inverse to function $h$. Using exponent -1 in such a manner makes sense in the unique case, namely for the multiplicative inverse, or reciprocal, $a^{-1}$ to element $a$, e.g. by numbers or matrices. In practically every context with inversion, this operation can be confused with exponentiation. Using the form of exponentiation makes no sense by inversion. Further using the exponent gives a formula an additional higher level by notation. Therefore, expressions with inverses in indexes and exponents cannot often be properly represented. However, this common notation of the inverse is traditional and could be further used along with a more suitable notation.

[Gelimson 2012b] proposed a possible alternative notation of the inverse. For example, the inverse to function
$$y = h(x)$$
is now denoted by
$$x = \underline{h}(y)$$
using the underline with no necessity of revocation. Simultaneously adding the traditional notation via
$$x = \underline{h}(y) = h^{-1}(y)$$
has some obvious advantages, too. It gives an aid to remember the sense of this new notation via building an association with the common notation which helps here this new notation. On the other hand, the new notation helps here the common notation. The reason is that using $\underline{h}(y)$ brings here at least doubt in the exponentiation role of inversion.

Let the same one-argument and one-value real-number function
$$y = h(x)$$
be namely nonstrictly monotonic. Then there are at least two different values $x_1$ and $x_2$ of argument $x$ such that
$$x_1 < x_2$$
but
$$y_1 = h(x_1) = h(x_2) = y_2 .$$
Let us denote their common value via
$$y_{12} = y_1 = y_2 .$$
Then image $y_{12}$ has at least two distinct preimages $x_1$ and $x_2$. Determine and gather all the distinct preimages of the same image $y_{12}$. They build a set which may be denoted by
$$\underline{h}(y_{12}).$$
In particular, this set contains both $x_1$ and $x_2$. Further there exist both the greatest lower bound
$$\underline{h}_{\inf}(y_{12})$$
and the least upper bound
$$\underline{h}_{\sup}(y_{12})$$
of this set. It seems to be better to denote them via $|\underline{h}$ and $\underline{h}|$, respectively. Namely,
$$|\underline{h}(y_{12}) = \underline{h}_{\inf}(y_{12}),$$
$$\underline{h}|(y_{12}) = \underline{h}_{\sup}(y_{12}).$$

If the both bounds are truly taken, then they are the minimal and the maximal elements of this set, respectively, with natural notation

$$\underline{h}_{min}(y_{12})$$

and

$$\underline{h}_{max}(y_{12}).$$

Then

$$|\underline{h}(y_{12}) = \underline{h}_{min}(y_{12}) = \underline{h}_{inf}(y_{12}),$$
$$\underline{h}|(y_{12}) = \underline{h}_{max}(y_{12}) = \underline{h}_{sup}(y_{12}).$$

Analyze in such a manner every image

$$y = y_{12}$$

with at least two distinct preimages $x_1$ and $x_2$ .

Further for every image y with the only preimage x ,

$$|\underline{h}(y) = \underline{h}_{min}(y) = \underline{h}_{inf}(y) = x$$

and

$$\underline{h}|(y) = \underline{h}_{max}(y) = \underline{h}_{sup}(y) = x .$$

Therefore, there exist the infimum (or greatest lower bound) inverse function

$$x = |\underline{h}(y)$$

and the supremum (or greatest lower bound) inverse function

$$x = \underline{h}|(y).$$

They are one-argument and one-value real-number functions and namely the both extreme functions among all the generally many-valued functions inverse to generally nonstrictly monotonically increasing one-argument and one-value real-number function

$$y = h(x).$$

Then there is an interval, one of two half-segments, or a segment

$$(|x , x|), (|x , x|], [|x , x|), [|x , x|]$$

whose either excluded or included (independently from one another) endpoints are

$$|x = |\underline{h}(y)$$

and

$$x| = \underline{h}|(y)$$

and which may be regarded as the total preimage of image y .

## 1.7. Integral (Cumulative) Probability Distribution Function Inversion

Let domain D of integral (cumulative) probability distribution function F(x) with range Ra = [0, 1] be one-dimensional Euclidean space R = (-∞, ∞). Then F(x) is a one-argument and one-value real-number function. It can be not only strictly monotonically increasing, but also locally nonstrictly monotonically increasing. If its arbitrary image $y_1$ belonging to range Ra = [0, 1] has at least two distinct preimages $x_1$ and $x_2$, then there is an interval, one of two half-segments, or a segment
$$(|\underline{F}(y_{12}), \underline{F}|(y_{12})), (|\underline{F}(y_{12}), \underline{F}|(y_{12})], [|\underline{F}(y_{12}), \underline{F}|(y_{12})), [|\underline{F}(y_{12}), \underline{F}|(y_{12})]$$
whose either excluded or included (independently from one another) endpoints are
$$x = |\underline{F}(y_{12})$$
and
$$x = \underline{F}|(y_{12})$$
and which may be regarded as the total preimage of image $y_{12}$.
Nota bene: On this interval possibly excluding its subset of zero measure, probability density function f(x) vanishes. Otherwise, integral (cumulative) probability distribution function F(x) could not be constant $y_{12}$.
Therefore, there exist the infimum inverse function
$$x = |\underline{F}(y)$$
and the supremum inverse function
$$x = \underline{F}|(y).$$
They are one-argument and one-value real-number functions and namely the both extreme functions among all the generally many-valued functions inverse to generally nonstrictly monotonically increasing one-argument and one-value real-number function
$$y = F(x).$$

## 1.8. Generally Nonstrictly Monotonic Sequence Continualization and Inversion

Let domain D of function g(x) be one-dimensional Euclidean space:
$$D = R = (-\infty, \infty).$$
Let namely finite real-number segment
$$S' = [a, b] \subset R = (-\infty, \infty) \ (-\infty < a < b < \infty)$$
be the extended support of function g(x) so that
$$g(x) = 0$$
at any
$$x \in R \setminus S' = (-\infty, a) \cup (b, \infty).$$
Let further n ($n \in N = \{1, 2, ...\}$) intermediate points $c_1, c_2, c_3, c_4, ..., c_{n-3}, c_{n-2}, c_{n-1}, c_n$ in the non-decreasing order so that
$$a \leq c_1 \leq c_2 \leq c_3 \leq c_4 \leq ... \leq c_{n-3} \leq c_{n-2} \leq c_{n-1} \leq c_n \leq b$$
divide this segment into $n + 1$ parts (pieces) of generally different lengths. To unify the notation, denote
$$c_0 = a,$$
$$c_{n+1} = b,$$
$$c(i) = c_i \ (i = 0, 1, 2, ..., n + 1).$$
Let us further generalize the last notation, namely for a generally nonstrictly monotonic sequence as function c(z) of index z
$$c(z) = c_z \ (z = z', z' + 1, z' + 2, z' + 3, z' + 4, z' + 5, ..., z'', z'' + 1, ..., z''')$$
(Fig. 1).

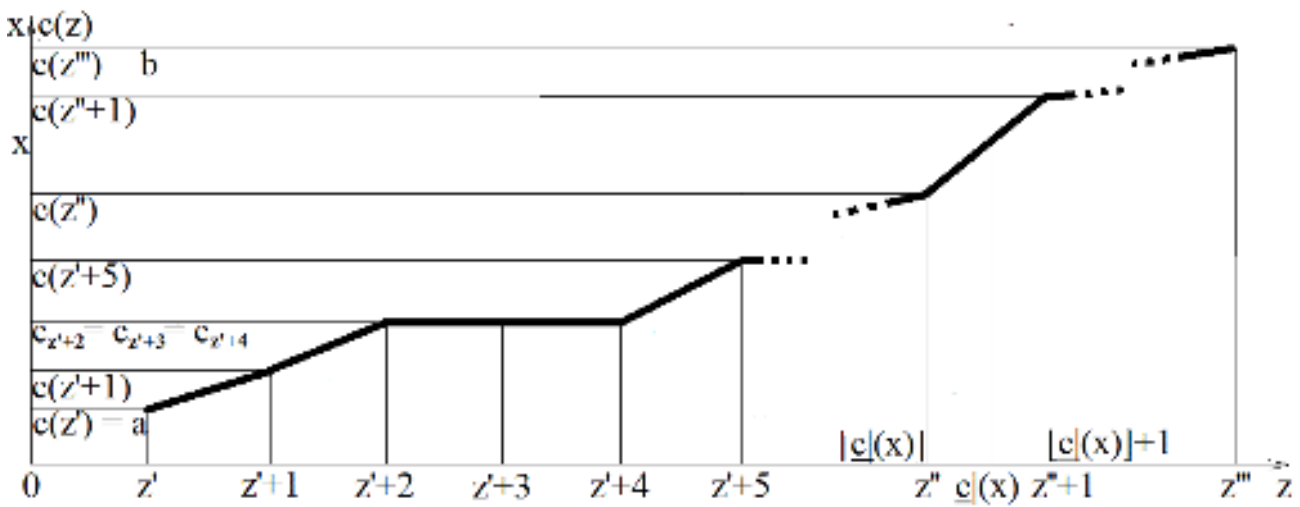

Figure 1. Generally nonstrictly monotonic sequence continualization and inversion

The numeration begins at any initial integer
$$z' \in Z = \{..., -2, -1, 0, 1, 2, ...\}$$
and ends at any final greater integer
$$z''' \in Z = \{..., -2, -1, 0, 1, 2, ...\},$$
$$z' < z''.$$
In the above particular case
$$z' = 0,$$
$$z''' = n + 1.$$

In our general case
$$a = c_{z'} \le c_{z'+1} \le c_{z'+2} \le c_{z'+3} \le c_{z'+4} \le c_{z'+5} \le \ldots \le c_{z''} \le c_{z''+1} \le \ldots \le c_{z'''} = b.$$
Nota bene: By namely strictly monotonically increasing sequence
$$c(z) = c_z,$$
all these inequalities would be strict:
$$a = c_{z'} < c_{z'+1} < c_{z'+2} < c_{z'+3} < c_{z'+4} < c_{z'+5} < \ldots < c_{z''} < c_{z''+1} < \ldots < c_{z'''} = b.$$
In this case only, all the $z''' - z' + 1$ points
$$a = c_{z'}, c_{z'+1}, c_{z'+2}, c_{z'+3}, c_{z'+4}, c_{z'+5}, \ldots, c_{z''}, c_{z''+1}, \ldots, c_{z'''} = b$$
would be namely distinct. However, in our general case, some of these $z''' - z' + 1$ points may coincide. In Figure 1, e.g., at least
$$c_{z'+2} = c_{z'+3} = c_{z'+4}$$
and
$$c_{z'+5} = c_{z'+6}$$
even if index $z' + 6$ is not explicitly shown but the direction from shown point $c_{z'+5}$ to omitted point $c_{z'+6}$ is obviously horizontal.

We shall use inverting this generally nonstrictly monotonic sequence as function $c(z)$ of index $z$ to provide space-saving notation of integral (cumulative) probability distribution function $F(x)$ with domain
$$R = (-\infty, \infty)$$
(one-dimensional Euclidean space) and range
$$Ra = [0, 1].$$
Then $F(x)$ is a one-argument and one-value real-number function. It can be not only strictly monotonically increasing, but also locally nonstrictly monotonically increasing. Further there are similar problems with inverting integral (cumulative) probability distribution function $F(x)$, in particular with determining quantiles, e.g. medians, quartiles, deciles, and percentiles. All this is our general aim. And our particular task here is as follows. For any real
$$x \in S' = [a, b] \subset R = (-\infty, \infty) \; (-\infty < a < b < \infty),$$
explicitly (closed-form) determine index $z''$ of left-closed and right-open interval
$$[c_{z''}, c_{z''+1}) = [c(z''), c(z'' + 1))$$
containing this x so that
$$x \in [c_{z''}, c_{z''+1}) = [c(z''), c(z'' + 1)),$$
$$c_{z''} \le x < c_{z''+1},$$
$$c(z'') \le x < c(z'' + 1).$$
Notata bene:
1. If the set of all the points $c(z)$ for which
$$c(z) = x$$
contains the only point, then denote this index $z$ via $z''$ :
$$z'' = \{z \in \{z', z' + 1, z' + 2, z' + 3, z' + 5, z' + 5, \ldots, z'''\} \mid c(z) = x\}.$$
As applied to such point $x$, this generally nonstrictly monotonic sequence as function $c(z)$ of index $z$ behaves so as if this sequence were strictly monotonic.
2. If the set of all the points $c(z)$ for which
$$c(z) = x$$
contains more than one point, then select namely the point with the maximal index $z$ among the indexes of all these points $c(z)$ and denote this maximal index $z$ via $z''$ :
$$z'' = \max\{z \in \{z', z' + 1, z' + 2, z' + 3, z' + 5, z' + 5, \ldots, z'''\} \mid c(z) = x\}.$$
3. If for real
$$x \in S' = [a, b] \subset R = (-\infty, \infty) \; (-\infty < a < b < \infty),$$
there is no point $c(z)$ for which
$$c(z) = x,$$
then it is also possible to explicitly (closed-form) determine index $z''$ of left-closed and right-open

interval
$$[c_{z''}, c_{z''+1}) = [c(z''), c(z''+1))$$
containing this x so that
$$x \in [c_{z''}, c_{z''+1}) = [c(z''), c(z''+1)),$$
$$c_{z''} \leq x < c_{z''+1},$$
$$c(z'') \leq x < c(z''+1).$$
Piecewise linear continualize sequence c(z) via simply consequently connecting all the z''' - z' + 1 diagram points in Figure 1
$$(z', c_{z'}), (z'+1, c_{z'+1}), (z'+2, c_{z'+2}), (z'+3, c_{z'+3}), (z'+4, c_{z'+4}), (z'+5, c_{z'+5}), \ldots, (z''' c_{z'''}).$$
Namely, on every real-number segment
$$[z, z+1]$$
for any integer
$$z \in \{z', z'+1, z'+2, z'+3, z'+5, z'+5, \ldots, z'''-1\},$$
apply linear interpolation
$$d(t) = c(z)(z+1-t) + c(z+1)(t-z)$$
$$= c_z(z+1-t) + c_{z+1}(t-z).$$
Notata bene:
1. For sequence c(z) which is an integer-argument real-valued function, its discrete domain is
$$\{z', z'+1, z'+2, z'+3, z'+5, z'+5, \ldots, z'''\}$$
containing these z''' - z' + 1 integer points only.
2. For real-argument real-valued function d(t), its continual domain is real segment
$$[z', z''']$$
also containing all the intermediate (internal) real points.
3. At every
$$z \in \{z', z'+1, z'+2, z'+3, z'+5, z'+5, \ldots, z'''\}$$
of these z''' - z' + 1 integer points, real-argument real-valued function d(t) coincides with sequence c(z) which is an integer-argument real-valued function.
Therefore, real-argument real-valued function d(t) is a generalization of sequence c(z) which is an integer-argument real-valued function with extending it from discrete domain
$$\{z', z'+1, z'+2, z'+3, z'+5, z'+5, \ldots, z'''\}$$
to continual domain
$$[z', z''']$$
also containing all the intermediate (internal) real points.
4. We may piecewise continualize sequence c(z) not only linearly but also nonlinearly with always piecewise conserving the monotonicity properties of sequence c(z) on every real-number segment
$$[z, z+1]$$
for any integer
$$z \in \{z', z'+1, z'+2, z'+3, z'+5, z'+5, \ldots, z'''-1\}.$$
5. However, if and only if
$$c_z = c_{z+1},$$
so that these points namely coincide, then the only possibility is to apply linear interpolation
$$d(t) = c_z(z+1-t) + c_{z+1}(t-z) =$$
$$= c_z(z+1-t) + c_z(t-z) = c_z = c_{z+1}$$
only (here even constant).
6. Otherwise, i.e. if and only if
$$c_z < c_{z+1},$$
so that these points are namely distinct, then along with linear interpolation
$$d(t) = c(z)(z+1-t) + c(z+1)(t-z)$$
$$= c_z(z+1-t) + c_{z+1}(t-z),$$
there are infinitely many different possibilities to apply nonlinear interpolation.

7. Moreover, due to nonlinear interpolation, there are infinitely many different possibilities to provide not only the continuity (which is the case by linear interpolation) of continualized function d(t), but also its differentiability.

8. To provide the differentiability of continualized function d(t) at such integer point
$$z \in \{z' + 1, z' + 2, z' + 3, z' + 5, z' + 5, \ldots, z''' - 1\}$$
that
$$c_{z-1} = c_z < c_{z+1},$$
it is necessary that the derivative of continualized function d(t) at this integer point from the right vanishes:
$$d'(z + 0) = d'(z - 0) = 0.$$

9. To provide the differentiability of continualized function d(t) at such integer point
$$z \in \{z' + 1, z' + 2, z' + 3, z' + 5, z' + 5, \ldots, z''' - 1\}$$
that
$$c_{z-1} < c_z = c_{z+1},$$
it is necessary that the derivative of continualized function d(t) at this integer point from the left vanishes:
$$d'(z - 0) = d'(z + 0) = 0.$$

10. To provide the differentiability of continualized function d(t) at such integer point
$$z \in \{z' + 1, z' + 2, z' + 3, z' + 5, z' + 5, \ldots, z''' - 1\}$$
that
$$c_{z-1} < c_z < c_{z+1},$$
i.e. z is a point of the two-sided strict increase of continualized function d(t), it is necessary that the derivatives of continualized function d(t) at this integer point from the left and from the right both have a common nonnegative value:
$$d'(z - 0) = d'(z + 0) \geq 0.$$

11. To provide the differentiability of continualized function d(t) at every integer point
$$z \in \{z', z' + 1, z' + 2, z' + 3, z' + 5, z' + 5, \ldots, z'''\},$$
it is sufficient that the derivative of continualized function d(t) at this integer point both from the left and from the right vanishes:
$$d'(z) = d'(z - 0) = d'(z + 0) = 0.$$

12. If
$$c_z < c_{z+1},$$
then there are infinitely many different possibilities to apply such nonlinear interpolation, e.g.
$$d(t) = (c_{z+1} - c_z)/2 \; |\sin[\pi(t - z - 1/2)]|^u \; \text{sign}(t - z - 1/2) + (c_z + c_{z+1})/2$$
for real segment [z, z + 1] by any integer
$$z \in \{z', z' + 1, z' + 2, z' + 3, z' + 5, z' + 5, \ldots, z''' - 1\}$$
and for any real
$$u \geq 1.$$

13. This curve
$$\{(t, d(t)) \mid t \in [z, z + 1]\}$$
is central symmetric about its inflection point
$$(z + 1/2, (c_z + c_{z+1})/2).$$
which is generally not necessary.

14. The inflection point of such a curve may be any in open rectangular
$$(z, z + 1) \times (c_z, c_{z+1})$$
without any symmetricity.

15. If
$$c_z < c_{z+1},$$
then there are infinitely many different possibilities to apply such nonlinear interpolation, e.g. via arbitrarily dividing each of the both pieces

$$[z, z+1]$$

and

$$[c_z, c_{z+1}]$$

into two subpieces

$$[z, z+q], [z+q, z+1]$$

and

$$[c_z, c_z + r(c_{z+1} - c_z)], [c_z + r(c_{z+1} - c_z), c_{z+1}].$$

This is equivalent to selecting any real

$$q \mid 0 < q < 1$$

and

$$r \mid 0 < r < 1.$$

Further define

$$d(t) = c_z + r(c_{z+1} - c_z)[(t - z)/q]^u, t \in [z, z+q],$$
$$d(t) = c_{z+1} - (1 - r)(c_{z+1} - c_z)[(z + 1 - t)/(1 - q)]^v, t \in [z+q, z+1]$$

for real segment $[z, z+1]$ by any integer

$$z \in \{z', z'+1, z'+2, z'+3, z'+5, z'+5, \ldots, z''' - 1\}.$$

and for any real

$$u \geq 1.$$

Hence we may rename function d to c and argument t to z and simply deal with real-argument real-valued function c(z).

Now we can solve our problem when for real

$$x \in S' = [a, b] \subset R = (-\infty, \infty) \ (-\infty < a < b < \infty),$$

there is no point c(z) for which

$$c(z) = x.$$

Apply inverting real-argument real-valued function c(z) to explicitly (closed-form) determine index z" of left-closed and right-open interval

$$[c_{z''}, c_{z''+1}) = [c(z''), c(z''+1))$$

containing this x so that

$$x \in [c_{z''}, c_{z''+1}) = [c(z''), c(z''+1)),$$
$$c_{z''} \leq x < c_{z''+1},$$
$$c(z'') \leq x < c(z''+1).$$

Denote the inverse function to real-argument real-valued function

$$x = c(z)$$

via

$$z = \underline{c}(x).$$

1. If the set of all the points c(z) for which

$$c(z) = x$$

contains the only point, then inverse function

$$z = \underline{c}(x)$$

provides the desired and required abscissa z.

2. If the set of all the points c(z) for which

$$c(z) = x$$

contains more than one point, then select namely the point with the maximal index z among the indexes of all these points c(z). To provide this, take namely

$$z = \underline{c}|(x).$$

However, we need namely index z" of left-closed and right-open interval

$$[c_{z''}, c_{z''+1}) = [c(z''), c(z''+1))$$

containing this x so that

$$x \in [c_{z''}, c_{z''+1}) = [c(z''), c(z''+1)),$$
$$c_{z''} \leq x < c_{z''+1},$$

$$c(z") \leq x < c(z"+1).$$

Using namely the supremum inverse function

$$z = \underline{c}|(x),$$

simply determine

$$z" = \lfloor \underline{c}|(x) \rfloor.$$

Here the floor (or entier) function [Encyclopaedia of Mathematics]

$$v = \lfloor w \rfloor$$

of real argument w gives the largest integer v less than or equal to w .

# 2. Piecewise Linear Probability Distribution

## 2.1. Main Definitions

Consider a general one-dimensional piecewise linear probability distribution (Fig. 2).

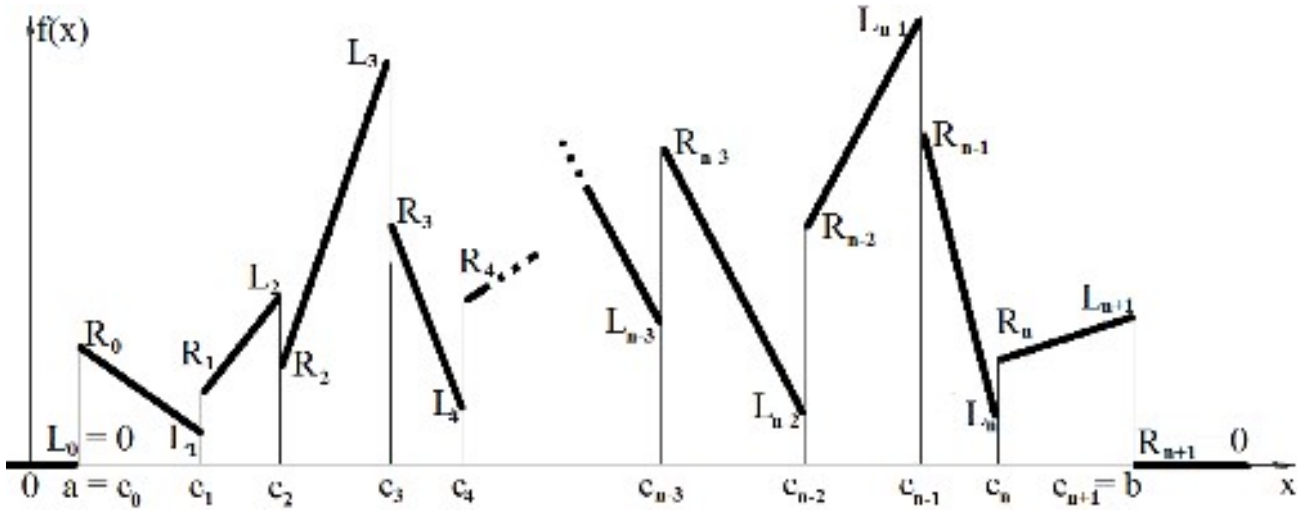

Figure 2. General one-dimensional piecewise linear probability distribution

Here probability density distribution function f(x) is as always non-negative everywhere ($-\infty < x < +\infty$) and can be positive on some so-called support which is a finite segment (closed interval)
$$-\infty < a \leq x \leq b < +\infty \ (a < b)$$
only. Let n ($n \in N = \{1, 2, ...\}$) intermediate points $c_1$, $c_2$, $c_3$, $c_4$, ... , $c_{n-3}$, $c_{n-2}$, $c_{n-1}$, $c_n$ in the non-decreasing order so that
$$a \leq c_1 \leq c_2 \leq c_3 \leq c_4 \leq ... \leq c_{n-3} \leq c_{n-2} \leq c_{n-1} \leq c_n \leq b$$
divide this segment into n + 1 parts (pieces) of generally different lengths. To unify the notation, denote
$$c_0 = a,$$
$$c_{n+1} = b,$$
$$c(i) = c_i \ (i = 0, 1, 2, ... , n + 1).$$
On each of n + 1 open intervals
$$c_i < x < c_{i+1} \ (i = 0, 1, 2, ... , n),$$
probability density distribution function f(x) is linear. At n + 2 points
$$c_i \ (i = 0, 1, 2, ... , n + 1),$$
f(x) may take any finite real values. The following considerations (possibly excepting mode values below) do not depend on these values. At each of n + 2 points
$$c_i \ (i = 0, 1, 2, ... , n + 1),$$
left and right one-sided limits
$$\lim f(x) = L_i \ (x \to c_i - 0),$$
$$\lim f(x) = R_i \ (x \to c_i + 0)$$
are any generally different finite real values. Naturally, we have
$$L_0 = 0,$$
$$R_{n+1} = 0.$$

Then on each of n + 1 open intervals
$$c_i < x < c_{i+1} \ (i = 0, 1, 2, \ldots, n),$$
linear probability density distribution function
$$f(x) = R_i + (L_{i+1} - R_i)(x - c_i)/(c_{i+1} - c_i)$$
$$= [R_i(c_{i+1} - x) + L_{i+1}(x - c_i)]/(c_{i+1} - c_i).$$
Using space-saving notation [Gelimson 2012a], represent non-negative-valued function f(x) via
$$f_{[0,\infty)}(x_{(-\infty,\infty)}) = 0_{(-\infty,c(0)) \cup (c(n+1),\infty)} \cup \cup_{i=0}^{n}[R_i(c_{i+1} - x) + L_{i+1}(x - c_i)]/(c_{i+1} - c_i)_{(c(i),c(i+1))} \cup \cup_{i=0}^{n+1} f_{\{c(i)\}}(c_i),$$
or, simplifying $f_{\{c(i)\}}(c_i)$ via identifying [Gelimson 2003a, 2003b] one-point set $\{c(i)\} = \{c_i\}$ at least here with this point $c(i) = c_i$ itself, via
$$f_{[0,\infty)}(x_{(-\infty,\infty)}) = 0_{(-\infty,c(0)) \cup (c(n+1),\infty)} \cup \cup_{i=0}^{n}[R_i(c_{i+1} - x) + L_{i+1}(x - c_i)]/(c_{i+1} - c_i)_{(c(i),c(i+1))} \cup \cup_{i=0}^{n+1} f_{c(i)}(c_i),$$
or, further simplifying $f_{c(i)}(c_i)$ via omitting index $c(i) = c_i$ coinciding with argument $c_i$, which is admissible if and only if argument $c_i$ is explicitly indicated, e.g. here (but NOT after replacing expression $f(c_i)$ via its value, e.g. a number), via
$$f_{[0,\infty)}(x_{(-\infty,\infty)}) = 0_{(-\infty,c(0)) \cup (c(n+1),\infty)} \cup \cup_{i=0}^{n}[R_i(c_{i+1} - x) + L_{i+1}(x - c_i)]/(c_{i+1} - c_i)_{(c(i),c(i+1))} \cup \cup_{i=0}^{n+1} f(c_i)$$
on the whole real axis $(-\infty, \infty)$
where
index $[0, \infty)$ in $f_{[0,\infty)}$ indicates the domain of dependent variable f and hence the range of function f(x),
index $(-\infty, \infty)$ in $x_{(-\infty,\infty)}$ indicates the range of x and hence the domain of one-argument function f(x),
index $(-\infty, c(0)) \cup (c(n+1), \infty) = (-\infty, c_0) \cup (c_{n+1}, \infty)$ in $0_{(-\infty,c(0))}$ indicates that function f(x) = 0 on its subdomain $(-\infty, c(0)) \cup (c(n+1), \infty) = (-\infty, c_0) \cup (c_{n+1}, \infty)$,
symbol $\cup$ unifies subfunctions on subdomains similarly to symbol $\cup$ in set theory and can be also indexed with an index range,
bounds 0 and n of index i in $\cup_{i=0}^{n}$ indicate that the range of index i is $\{0, 1, 2, \ldots, n\}$,
index $(c(i), c(i+1)) = (c_i, c_{i+1})$ in $\{[R_i(c_{i+1} - x) + L_{i+1}(x - c_i)]/(c_{i+1} - c_i)\}_{(c(i), c(i+1))}$ indicates that function $f(x) = [R_i(c_{i+1} - x) + L_{i+1}(x - c_i)]/(c_{i+1} - c_i)$ on its subdomain $(c(i), c(i+1)) = (c_i, c_{i+1})$,
index $\{c(i)\} = \{c_i\}$ in $f_{\{c(i)\}}(c_i)$ indicates that function $f(x) = f(c_i)$ on its subdomain $\{c(i)\} = \{c_i\}$.

## 2.2. Normalization Condition

The probability of the event that X takes any finite real value is namely 1 because this event is certain. This gives integral normalization condition
$$\int_{-\infty}^{+\infty} f(x)dx = 1.$$
In our case we have
$$1 = \int_{-\infty}^{+\infty} f(x)dx = \int_a^b f(x)dx$$
$$= \Sigma_{i=0}^n \int_{c_{(i)}}^{c_{(i+1)}} f(x)dx = \Sigma_{i=0}^n \int_{c_{(i)}}^{c_{(i+1)}} [R_i(c_{i+1} - x) + L_{i+1}(x - c_i)]/(c_{i+1} - c_i) \, dx$$
$$= \Sigma_{i=0}^n \{R_i[c_{i+1}(c_{i+1} - c_i) - (c_{i+1}^2 - c_i^2)/2] + L_{i+1}[(c_{i+1}^2 - c_i^2)/2 - c_i(c_{i+1} - c_i)]\}/(c_{i+1} - c_i)$$
$$= \Sigma_{i=0}^n \{R_i[c_{i+1} - (c_{i+1} + c_i)/2] + L_{i+1}[(c_{i+1} + c_i)/2 - c_i]\}$$
$$= \Sigma_{i=0}^n (R_i + L_{i+1})(c_{i+1} - c_i)/2.$$
We can also obtain this result at once rather geometrically than analytically, namely via adding the areas of the n + 1 rectangular trapezoids.
Therefore, to provide a possible (an admissible) probability density distribution function, necessary and sufficient integral normalization condition
$$\Sigma_{i=0}^n (R_i + L_{i+1})(c_{i+1} - c_i) = 2$$
has to be satisfied.

## 2.3. Normalization Algorithm

Nota bene: The obtained normalization condition is one condition only for
$$(n + 1) + (n + 1) + (n + 2) = 3n + 4$$
unknowns
$$R_i \ (i = 0, 1, 2, \ldots, n),$$
$$L_i \ (i = 1, 2, 3, \ldots, n + 1),$$
$$c_i \ (i = 0, 1, 2, \ldots, n + 1).$$
Additionally,
$$R_i \geq 0 \ (i = 0, 1, 2, \ldots, n),$$
$$L_i \geq 0 \ (i = 1, 2, 3, \ldots, n + 1),$$
$$c_0 \leq c_1 \leq c_2 \leq c_3 \leq c_4 \leq \ldots \leq c_{n-3} \leq c_{n-2} \leq c_{n-1} \leq c_n \leq c_{n+1}.$$
Generally, it is not possible to simply take any admissible values of
$$3n + 4 - 1 = 3n + 3$$
unknowns and then to determine the value of the remaining unknown via this condition because it can happen that this value is inadmissible.

A natural idea, way, and algorithm to avoid this difficulty are as follows:

1. Fix
$$c_0 \leq c_1 \leq c_2 \leq c_3 \leq c_4 \leq \ldots \leq c_{n-3} \leq c_{n-2} \leq c_{n-1} \leq c_n \leq c_{n+1}.$$

2. Take any
$$R_i' \geq 0 \ (i = 0, 1, 2, \ldots, n),$$
$$L_i' \geq 0 \ (i = 1, 2, 3, \ldots, n + 1)$$
so that there is at least one namely positive number among these $2n + 2$ non-negative numbers.

3. Let
$$R_i \ (i = 0, 1, 2, \ldots, n),$$
$$L_i \ (i = 1, 2, 3, \ldots, n + 1)$$
be proportional to
$$R_i' \geq 0 \ (i = 0, 1, 2, \ldots, n),$$
$$L_i' \geq 0 \ (i = 1, 2, 3, \ldots, n + 1),$$
respectively, with a common namely positive factor k so that
$$R_i = kR_i' \ (i = 0, 1, 2, \ldots, n),$$
$$L_i = kL_i' \ (i = 1, 2, 3, \ldots, n + 1).$$

4. Explicitly determine the value of parameter k as the only unknown via this necessary and sufficient integral normalization condition
$$\Sigma_{i=0}^n (R_i + L_{i+1})(c_{i+1} - c_i) = 2$$
so that
$$k = 2 / \Sigma_{i=0}^n (R_i' + L_{i+1}')(c_{i+1} - c_i).$$

5. Explicitly determine
$$R_i = kR_i' \ (i = 0, 1, 2, \ldots, n),$$
$$L_i = kL_i' \ (i = 1, 2, 3, \ldots, n + 1).$$

## 2.4. Integral (Cumulative) Probability Distribution Function

Integral (cumulative) probability distribution function
$$F(x) = P(X \leq x) = \int_{-\infty}^{x} f(t)dt$$
is probability $P(X \leq x)$ that real-number random variable X takes a real-number value not greater than x .
For $x \leq a = c_0$ , this definition gives
$$F(x) = P(X \leq x) = \int_{-\infty}^{x} f(t)dt = 0.$$
For $x \geq b = c_{n+1}$ , this definition gives
$$F(x) = P(X \leq x) = \int_{-\infty}^{x} f(t)dt = 1.$$
For $c_0 = a < x < b = c_{n+1}$ , we can use the following natural idea, way, and algorithm:
1. Determine such value j of index i that
$$c_j \leq x < c_{j+1} .$$
There exists such value j and namely the only. Indeed, consider set
$$\{i \mid i \in \{0, 1, 2, ... , n\}, c_i \leq x\}.$$
It is nonempty because it contains at least 0 for which
$$c_0 = a < x$$
and hence
$$c_0 \leq x .$$
It is finite and strictly ordered by relation < . Therefore, there exists its maximal element
$$j = \max\{i \mid i \in \{0, 1, 2, ... , n\}, c_i \leq x\},$$
and this maximal element is namely the only. Then for this maximal element j , in addition to
$$c_j \leq x ,$$
inequality
$$x < c_{j+1}$$
also holds. Indeed, otherwise, we would have
$$x \geq c_{j+1}$$
and
$$c_{j+1} \leq x ,$$
so that this j could not be namely the maximal element of this set.
2. Determine
$$F(x) = P(X \leq x) = \int_{-\infty}^{x} f(t)dt = \int_{a}^{x} f(t)dt = \int_{c(0)}^{x} f(t)dt$$
$$= \Sigma_{i=0}^{j-1} \int_{c(i)}^{c(i+1)} f(t)dt + \int_{c(j)}^{x} f(t)dt$$
$$= \Sigma_{i=0}^{j-1} \int_{c(i)}^{c(i+1)} [R_i(c_{i+1} - t) + L_{i+1}(t - c_i)]/(c_{i+1} - c_i) \, dt$$
$$+ \int_{c(j)}^{x} [R_j(c_{j+1} - t) + L_{j+1}(t - c_j)]/(c_{j+1} - c_j) \, dt$$
$$= \Sigma_{i=0}^{j-1} \{R_i[c_{i+1}(c_{i+1} - c_i) - (c_{i+1}^2 - c_i^2)/2] + L_{i+1}[(c_{i+1}^2 - c_i^2)/2 - c_i(c_{i+1} - c_i)]\}/(c_{i+1} - c_i)$$
$$+ \{R_j[c_{j+1}(x - c_j) - (x^2 - c_j^2)/2] + L_{j+1}[(x^2 - c_j^2)/2 - c_j(x - c_j)]\}/(c_{j+1} - c_j)$$
$$= \Sigma_{i=0}^{j-1} \{R_i[c_{i+1} - (c_{i+1} + c_i)/2] + L_{i+1}[(c_{i+1} + c_i)/2 - c_i]\}$$
$$+ \{R_j[c_{j+1}(x - c_j) - (x^2 - c_j^2)/2] + L_{j+1}[(x^2 - c_j^2)/2 - c_j(x - c_j)]\}/(c_{j+1} - c_j)$$
$$= \Sigma_{i=0}^{j-1} (R_i + L_{i+1})(c_{i+1} - c_i)/2 + \{Rj[c_{j+1}(x - c_j) - (x^2 - c_j^2)/2] + L_{j+1}[(x^2 - c_j^2)/2 - c_j(x - c_j)]\}/(c_{j+1} - c_j).$$

## 1.4. Mean Value (Mathematical Expectation)

Use the common integral definition [Cramér] of the mean value (mathematical expectation)
$$\mu = E(X) = \int_{-\infty}^{+\infty} xf(x)dx .$$
In our case we determine
$$\mu = \int_{-\infty}^{+\infty} xf(x)dx = \int_a^b xf(x)dx$$
$$= \Sigma_{i=0}^n \int_{c_{(i)}}^{c_{(i+1)}} xf(x)dx = \Sigma_{i=0}^n \int_{c_{(i)}}^{c_{(i+1)}} [R_i(c_{i+1}x - x^2) + L_{i+1}(x^2 - c_ix)]/(c_{i+1} - c_i)\, dx$$
$$= \Sigma_{i=0}^n \{R_i[c_{i+1}(c_{i+1}^2 - c_i^2)/2 - (c_{i+1}^3 - c_i^3)/3] + L_{i+1}[(c_{i+1}^3 - c_i^3)/3 - c_i(c_{i+1}^2 - c_i^2)/2]/(c_{i+1} - c_i)\}$$
$$= 1/6\, \Sigma_{i=0}^n \{R_i[3c_{i+1}(c_{i+1} + c_i) - 2(c_{i+1}^2 + c_{i+1}c_i + c_i^2)] + L_{i+1}[2(c_{i+1}^2 + c_{i+1}c_i + c_i^2) - 3c_i(c_{i+1} + c_i)]\}$$
$$= 1/6\, \Sigma_{i=0}^n [R_i(c_{i+1}^2 + c_{i+1}c_i - 2c_i^2) + L_{i+1}(2c_{i+1}^2 - c_{i+1}c_i - c_i^2)]$$
$$= 1/6\, \Sigma_{i=0}^n (c_{i+1} - c_i)[R_i(c_{i+1} + 2c_i) + L_{i+1}(2c_{i+1} + c_i)]$$
and, finally,
$$\mu = \Sigma_{i=0}^n (c_{i+1} - c_i)[R_i(2c_i + c_{i+1}) + L_{i+1}(c_i + 2c_{i+1})]/6.$$

## 2.5. Median Values

Use the common integral definition [Cramér] of median values $v$ for any of which both
$$P(X \leq v) \geq 1/2$$
and
$$P(X \geq v) \geq 1/2.$$
For a continual real-number random variable $X$,
$$P(X \leq v) = \int_{-\infty}^{v} f(x)dx = P(X \geq v) = \int_{v}^{+\infty} f(x)dx = 1/2.$$
To determine the set of all the median values $v$, we can use the following natural idea, way, and algorithm:
1. First consider
$$c_i \ (i = 0, 1, 2, \ldots, n+1)$$
not far from $\mu$ and determine both
$$L = \max\{i \mid \int_{-\infty}^{c(i)} f(x)dx < 1/2\}$$
and
$$R = \min\{i \mid \int_{c(i)}^{+\infty} f(x)dx < 1/2\}.$$
Then both
$$\int_{-\infty}^{c(L+1)} f(x)dx \geq 1/2$$
and
$$\int_{c(R-1)}^{+\infty} f(x)dx \geq 1/2.$$
2. On half-closed interval
$$c(L) = c_L < v \leq c_{L+1} = c(L+1),$$
determine
$$v_{min} = \inf\{v \mid \int_{-\infty}^{v} f(x)dx = 1/2\}.$$
3. On half-closed interval
$$c(R-1) = c_{R-1} \leq v < c_R = c(R),$$
determine
$$v_{max} = \sup\{v \mid \int_{v}^{+\infty} f(x)dx = 1/2\}.$$
4. Then the set of all the median values $v$ is the interval whose endpoints are
$$v_{min} \leq v_{max}$$
each of which is included into the interval if and only if the corresponding greatest lower and/or least upper bound is really taken so that
$$v_{min} = \min\{v \mid \int_{-\infty}^{v} f(x)dx = 1/2\}$$
and/or
$$v_{max} = \max\{v \mid \int_{v}^{+\infty} f(x)dx = 1/2\},$$
respectively.
Notata bene:
1. If
$$v_{min} = v_{max},$$
then the corresponding greatest lower and/or least upper bound is really taken so that
$$v_{min} = \min\{v \mid \int_{-\infty}^{v} f(x)dx = 1/2\}$$
and
$$v_{max} = \max\{v \mid \int_{v}^{+\infty} f(x)dx = 1/2\},$$
hence the closed interval
$$v_{min} \leq v \leq v_{max}$$
contains the only median value
$$v = v_{min} = v_{max}.$$
2. If
$$v_{min} < v_{max},$$

then the integral of f(x) on the interval whose endpoints are $v_{min}$ and $v_{max}$ vanishes independently of their including or excluding. Hence on this interval, non-negative probability density distribution function f(x) also vanishes possibly excepting points whose set has zero measure (in our case, a finite set).

## 2.6. Mode Values

To begin with, consider the common definition [Cramér] of mode values for any of which probability density distribution function f(x) takes its maximum value $f_{max}$. For continual distributions, generalize this definition in the following directions:
1. Replace the maximum value $f_{max}$ with the supremum value $f_{sup}$ which always exists. The reason is that it is possible (for piecewise linear probability distributions, too) that function f(x) is discontinuous and does not take the supremum value $f_{sup}$ so that the maximum value $f_{max}$ does not exist at all.
2. Extend the range of function f(x), i.e. the set of values function f(x) really (truly) takes, via all the limiting points of this set. Then the extended range is a closed set and contains, in particular, the supremum value $f_{sup}$.
3. Extend the domain of function f(x), i.e. the set of points at which function f(x) is properly defined, via all the limiting points of this set. Then the extended domain is a closed set which contains all its limiting points.
4. Admit modes to also correspond to the one-sided limits of function f(x) separately if necessary. This is important for discontinuous function f(x) with jumps.
5. At any interval endpoint $c_i$, along with the given value of $f(c_i)$, take into account the one-sided limits $L_i$ and $R_i$ of function f(x), e.g. any of the following reasonable options for value $f(c_i)$:
5.1. Take the given value of $f(c_i)$ itself.
5.2. Take
$$f(c_i) = \max\{L_i, R_i\}.$$
5.3. Take
$$f(c_i) = (L_i + R_i)/2.$$
6. At any interval endpoint $c_i$, along with $c_i$ itself, take into account the one-sided limiting points $c_i - 0$ and $c_i + 0$ corresponding to one-sided limits $L_i$ and $R_i$ of function f(x), respectively, e.g. any of the following reasonable options for $c_i$:
6.1. Take the given value of $c_i$ itself.
6.2. For modes, rather than $c_i$, consider
$$c_i - 0 \text{ if } L_i > R_i,$$
$$c_i + 0 \text{ if } L_i < R_i,$$
and quantiset [Gelimson 2003a, 2003b]
$$\{_{1/2}(c_i - 0), _{1/2}(c_i + 0)\}° \text{ if } L_i = R_i.$$
This quantiset consists of two quantielements
$$_{1/2}(c_i - 0), _{1/2}(c_i + 0)$$
with bases
$$c_i - 0, c_i + 0,$$
respectively.
Here each of elements $c_i - 0$ and $c_i + 0$ has quantity 1/2 so that the total unit quantity is equally divided between these both elements.
In particular, for a piecewise linear probability distribution with probability density function f(x), anyone of the following values can reasonably play the role of $f_{sup}$:
$$\max\{\max\{f(c_i) \mid i = 0, 1, 2, \ldots, n+1\}, \max\{L_i \mid i = 0, 1, \ldots, n+1\}, \max\{R_i \mid i = 0, 1, \ldots, n+1\}\},$$
$$\max\{\max\{f(c_i) \mid i = 0, 1, 2, \ldots, n+1\}, \max\{(L_i + R_i)/2 \mid i = 0, 1, 2, \ldots, n+1\}\},$$
$$\max\{\max\{L_i \mid i = 0, 1, 2, \ldots, n+1\}, \max\{R_i \mid i = 0, 1, \ldots, n+1\}\},$$
$$\max\{(L_i + R_i)/2 \mid i = 0, 1, 2, \ldots, n+1\}.$$
If $f(c_i) = f_{sup}$ at some i, then $c_i$ at this i is one of the modes.
If $L_i = f_{sup}$ at some i, then $c_i - 0$ at this i is one of the modes.
If $R_i = f_{sup}$ at some i, then $c_i + 0$ at this i is one of the modes.

If $(L_i + R_i)/2 = f_{sup}$ at some $i$, then quantiset
$$\{_{1/2}(c_i - 0), \,_{1/2}(c_i + 0)\}^\circ$$
at this $i$ is one of the modes.

Nota bene: The set of all the modes contains the corresponding separate points $c_i$, as well as one-sided limits $c_i - 0$ and $c_i + 0$, and includes open intervals
$$c_i < x < c_{i+1} \,(i = 1, 2, \ldots, n - 1)$$
for which
$$R_i = L_{i+1} = f_{sup} \,.$$

## 2.7. Variance

Use the common integral definition [Cramér] of the variance $\sigma^2$ of a random variable X as its second central moment, namely the squared standard deviation $\sigma$, or the expected value of the squared deviation from the mean:

$$\sigma^2 = E[(X - \mu)^2] = \int_{-\infty}^{+\infty} (x - \mu)^2 f(x) dx .$$

In our case we determine

$$\sigma^2 = \int_{-\infty}^{+\infty} (x - \mu)^2 f(x) dx = \int_a^b (x - \mu)^2 f(x) dx = \Sigma_{i=0}^n \int_{c_{(i)}}^{c_{(i+1)}} (x - \mu)^2 f(x) dx$$

$$= \Sigma_{i=0}^n \int_{c_{(i)}}^{c_{(i+1)}} [R_i(x^2 - 2\mu x + \mu^2)(c_{i+1} - x) + L_{i+1}(x^2 - 2\mu x + \mu^2)(x - c_i)]/(c_{i+1} - c_i) \, dx$$

$$= \Sigma_{i=0}^n \int_{c_{(i)}}^{c_{(i+1)}} \{R_i[-x^3 + (2\mu + c_{i+1})x^2 - (\mu^2 + 2\mu c_{i+1})x + \mu^2 c_{i+1}]$$
$$+ L_{i+1}[x^3 - (2\mu + c_i)x^2 + (\mu^2 + 2\mu c_i)x - \mu^2 c_i]\}/(c_{i+1} - c_i) \, dx$$

$$= \Sigma_{i=0}^n \{R_i[-(c_{i+1}^4 - c_i^4)/4 + (2\mu + c_{i+1})(c_{i+1}^3 - c_i^3)/3 - (\mu^2 + 2\mu c_{i+1})(c_{i+1}^2 - c_i^2)/2 + \mu^2 c_{i+1}(c_{i+1} - c_i)]$$
$$+ L_{i+1}[(c_{i+1}^4 - c_i^4)/4 - (2\mu + c_i)(c_{i+1}^3 - c_i^3)/3 + (\mu^2 + 2\mu c_i)(c_{i+1}^2 - c_i^2)/2 - \mu^2 c_i(c_{i+1} - c_i)]\}/(c_{i+1} - c_i)$$

$$= 1/12 \, \Sigma_{i=0}^n \{R_i[-3c_{i+1}^3 - 3c_{i+1}^2 c_i - 3c_{i+1} c_i^2 - 3c_i^3 + (4c_{i+1} + 8\mu)(c_{i+1}^2 + c_{i+1} c_i + c_i^2) - (12\mu c_{i+1} + 6\mu^2)(c_{i+1} + c_i) + 12\mu^2 c_{i+1}]$$
$$+ L_{i+1}[3c_{i+1}^3 + 3c_{i+1}^2 c_i + 3c_{i+1} c_i^2 + 3c_i^3 - (4c_i + 8\mu)(c_{iївs+1}^2 + c_{i+1} c_i + c_i^2) + (12\mu c_i + 6\mu^2)(c_{i+1} + c_i) - 12\mu^2 c_i]\}$$

$$= 1/12 \, \Sigma_{i=0}^n [R_i(-3c_{i+1}^3 - 3c_{i+1}^2 c_i - 3c_{i+1} c_i^2 - 3c_i^3 + 4c_{i+1}^3 + 4c_{i+1}^2 c_i + 4c_{i+1} c_i^2 + 8\mu c_{i+1}^2 + 8\mu c_{i+1} c_i + 8\mu c_i^2 - 12\mu c_{i+1}^2 - 12\mu c_{i+1} c_i - 6\mu^2 c_{i+1} - 6\mu^2 c_i + 12\mu^2 c_{i+1})$$
$$+ L_{i+1}(3c_{i+1}^3 + 3c_{i+1}^2 c_i + 3c_{i+1} c_i^2 + 3c_i^3 - 4c_{i+1}^2 c_i - 4c_{i+1} c_i^2 - 4c_i^3 - 8\mu c_{i+1}^2 - 8\mu c_{i+1} c_i - 8\mu c_i^2 + 12\mu c_{i+1} c_i + 12\mu c_i^2 + 6\mu^2 c_{i+1} + 6\mu^2 c_i - 12\mu^2 c_i)]$$

$$= 1/12 \, \Sigma_{i=0}^n [R_i(c_{i+1}^3 + c_{i+1}^2 c_i + c_{i+1} c_i^2 - 3c_i^3 - 4\mu c_{i+1}^2 - 4\mu c_{i+1} c_i + 8\mu c_i^2 + 6\mu^2 c_{i+1} - 6\mu^2 c_i)$$
$$+ L_{i+1}(3c_{i+1}^3 - c_{i+1}^2 c_i - c_{i+1} c_i^2 - c_i^3 - 8\mu c_{i+1}^2 + 4\mu c_{i+1} c_i + 4\mu c_i^2 + 6\mu^2 c_{i+1} - 6\mu^2 c_i)]$$

$$= 1/12 \, \Sigma_{i=0}^n (c_{i+1} - c_i)[R_i(c_{i+1}^2 + c_{i+1} c_i + c_i^2 + c_{i+1} c_i + c_i^2 + c_i^2 - 4\mu c_{i+1} - 4\mu c_i - 4\mu c_i + 6\mu^2)$$
$$+ L_{i+1}(c_{i+1}^2 + c_{i+1} c_i + c_i^2 + c_{i+1} c_i + c_{i+1}^2 + c_{i+1}^2 - 4\mu c_{i+1} - 4\mu c_{i+1} - 4\mu c_i + 6\mu^2)]$$

$$= 1/12 \, \Sigma_{i=0}^n (c_{i+1} - c_i)[R_i(c_{i+1}^2 + 2c_{i+1} c_i + 3c_i^2 - 4\mu c_{i+1} - 8\mu c_i + 6\mu^2)$$
$$+ L_{i+1}(3c_{i+1}^2 + 2c_{i+1} c_i + c_i^2 - 8\mu c_{i+1} - 4\mu c_i + 6\mu^2)]$$

and, finally,

$$\sigma^2 = \Sigma_{i=0}^n (c_{i+1} - c_i)[R_i(c_{i+1}^2 + 2c_{i+1} c_i + 3c_i^2 - 4\mu c_{i+1} - 8\mu c_i + 6\mu^2)$$
$$+ L_{i+1}(3c_{i+1}^2 + 2c_{i+1} c_i + c_i^2 - 8\mu c_{i+1} - 4\mu c_i + 6\mu^2)]/12$$

where

$$\mu = \Sigma_{i=0}^n (c_{i+1} - c_i)[R_i(2c_i + c_{i+1}) + L_{i+1}(c_i + 2c_{i+1})]/6.$$

Nota bene: Similarly, we can also determine further initial and central moments etc. [Cramér], e.g. skewness

$$\gamma_1 = E[(X - \mu)^3/\sigma^3]$$

and excess

$$\gamma_2 = E[(X - \mu)^4/\sigma^4] - 3.$$

# 3. General Polygonal, or Piecewise Linear Continuous, Probability Distribution

## 3.1. Main Definitions

Consider a general one-dimensional polygonal, or piecewise linear continuous, probability distribution (Fig. 2) as a particular case of a general one-dimensional piecewise linear probability distribution.

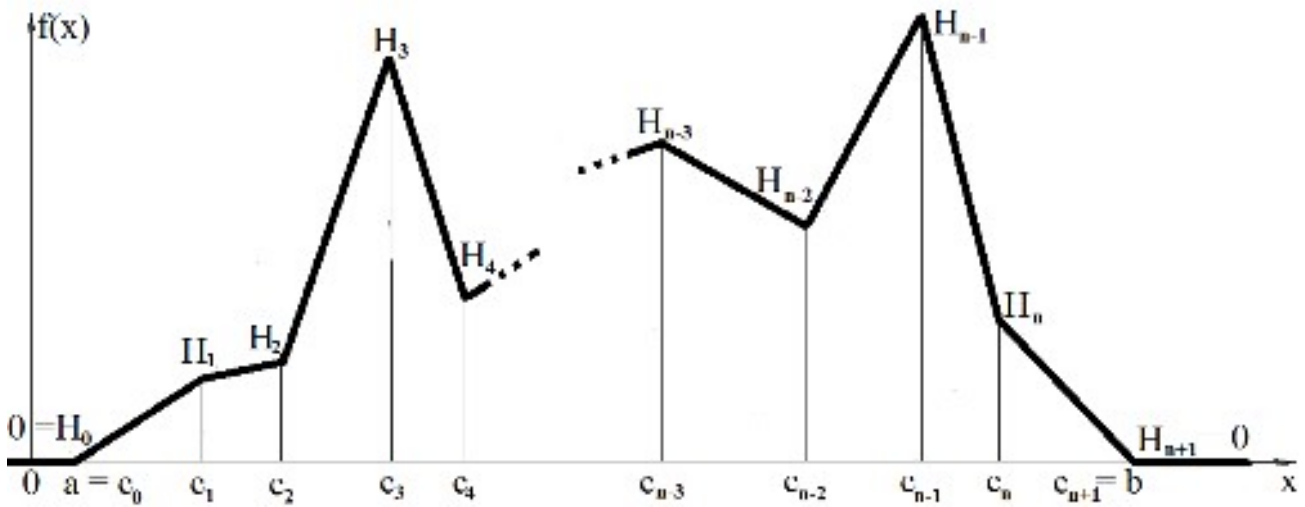

Figure 3. General one-dimensional piecewise linear continuous probability distribution

Here probability density distribution function $f(x)$ is as always non-negative everywhere ($-\infty < x < +\infty$) and can be positive on some finite segment (closed interval)
$$-\infty < a \leq x \leq b < +\infty \ (a < b)$$
only. Let $n$ ($n \in N = \{1, 2, ...\}$) intermediate points $c_1, c_2, c_3, c_4, \ldots, c_{n-3}, c_{n-2}, c_{n-1}, c_n$ in the non-decreasing order so that
$$a \leq c_1 \leq c_2 \leq c_3 \leq c_4 \leq \ldots \leq c_{n-3} \leq c_{n-2} \leq c_{n-1} \leq c_n \leq b$$
divide this segment into $n + 1$ parts (pieces) of generally different lengths. To unify the notation, denote
$$c_0 = a,$$
$$c_{n+1} = b,$$
$$c(i) = c_i \ (i = 0, 1, 2, \ldots, n + 1).$$
On each of $n + 1$ closed intervals
$$c_i \leq x \leq c_{i+1} \ (i = 0, 1, 2, \ldots, n),$$
probability density distribution function $f(x)$ is linear. At $n + 2$ points
$$c_i \ (i = 0, 1, 2, \ldots, n + 1),$$
$f(x)$ takes finite non-negative values
$$H_i = f(c_i),$$
respectively. Naturally, we have
$$H_0 = 0,$$
$$H_{n+1} = 0.$$
Note that

$$H_i = f(c_i) \ (i = 1, 2, \ldots, n)$$

may be any finite non-negative values. At each of $n + 2$ points

$$c_i \ (i = 0, 1, 2, \ldots, n + 1),$$

left and right one-sided limits

$$\lim f(x) = L_i \ (x \to c_i - 0),$$
$$\lim f(x) = R_i \ (x \to c_i + 0)$$

are equal to one another and coincide with $f(c_i)$. Therefore, we obtain

$$H_i = L_i = R_i \ (i = 0, 1, 2, \ldots, n + 1),$$

which makes it possible to apply the above formulas for a piecewise linear probability distribution to a piecewise linear continuous probability distribution.

Then on each of $n + 1$ closed intervals

$$c_i \leq x \leq c_{i+1} \ (i = 0, 1, 2, \ldots, n),$$

linear probability density distribution function

$$f(x) = H_i + (H_{i+1} - H_i)(x - c_i)/(c_{i+1} - c_i)$$
$$= [H_i(c_{i+1} - x) + H_{i+1}(x - c_i)]/(c_{i+1} - c_i).$$

Use space-saving notation [Gelimson 2013] and the corresponding formula for a piecewise linear probability distribution. Then in our continuous case we represent non-negative-valued function $f(x)$ via

$$f_{[0, \infty)}(x_{(-\infty, \infty)}) = 0_{(-\infty, c(0)) \cup [c(n+1), \infty)} \cup \cup_{i=0}^{n} [H_i(c_{i+1} - x) + H_{i+1}(x - c_i)]/(c_{i+1} - c_i)_{[c(i), c(i+1))}$$

on the whole real axis $(-\infty, \infty)$.

Integral (cumulative) probability distribution function

$$F(x) = P(X \leq x) = \int_{-\infty}^{x} f(t)dt$$

is probability $P(X \leq x)$ that real-number random variable $X$ takes a real-number value not greater than $x$.

## 3.2. Normalization Condition

The probability of the event that X takes any finite real value is namely 1 because this event is certain. This gives integral normalization condition
$$\int_{-\infty}^{+\infty} f(x)dx = 1.$$
Use the corresponding formula for a piecewise linear probability distribution. Then in our continuous case we determine
$$1 = \int_{-\infty}^{+\infty} f(x)dx = \int_a^b f(x)dx$$
$$= \Sigma_{i=0}^n (R_i + L_{i+1})(c_{i+1} - c_i)/2$$
$$= \Sigma_{i=0}^n (H_i + H_{i+1})(c_{i+1} - c_i)/2$$
$$= \Sigma_{i=0}^n H_i(c_{i+1} - c_i)/2 + \Sigma_{i=0}^n H_{i+1}(c_{i+1} - c_i)/2.$$
We can also obtain this result at once rather geometrically than analytically, namely via adding the areas of the n + 1 rectangular trapezoids, among them 2 rectangular triangles at the endpoints a and b .
Now use
$$H_0 = 0,$$
$$H_{n+1} = 0.$$
Then
$$1 = \int_{-\infty}^{+\infty} f(x)dx = \Sigma_{i=1}^n H_i(c_{i+1} - c_i)/2 + \Sigma_{i=1}^n H_i(c_i - c_{i-1})/2$$
$$= \Sigma_{i=1}^n H_i(c_{i+1} - c_{i-1})/2.$$
Therefore, to provide a possible (an admissible) probability density distribution function, necessary and sufficient integral normalization condition
$$\Sigma_{i=1}^n H_i(c_{i+1} - c_{i-1}) = 2$$
has to be satisfied.

# 3.3. Normalization Algorithm

Nota bene: The obtained normalization condition is one condition only for
$$n + (n + 2) = 2n + 2$$
unknowns
$$H_i \ (i = 1, 2, \ldots, n),$$
$$c_i \ (i = 0, 1, 2, \ldots, n + 1).$$
Additionally,
$$H_i \geq 0 \ (i = 1, 2, 3, \ldots, n),$$
$$c_0 \leq c_1 \leq c_2 \leq c_3 \leq c_4 \leq \ldots \leq c_{n-3} \leq c_{n-2} \leq c_{n-1} \leq c_n \leq c_{n+1}.$$
Generally, it is not possible to simply take any admissible values of
$$2n + 2 - 1 = 2n + 1$$
unknowns and then to determine the value of the remaining unknown via this condition because it can happen that this value is inadmissible.
A natural idea, way, and algorithm to avoid this difficulty are as follows:
1. Fix
$$c_0 \leq c_1 \leq c_2 \leq c_3 \leq c_4 \leq \ldots \leq c_{n-3} \leq c_{n-2} \leq c_{n-1} \leq c_n \leq c_{n+1}.$$
2. Take any
$$H_i' \geq 0 \ (i = 1, 2, \ldots, n)$$
so that there is at least one namely positive number among these n non-negative numbers.
3. Let
$$H_i \ (i = 1, 2, 3, \ldots, n)$$
be proportional to
$$H_i' \geq 0 \ (i = 1, 2, 3, \ldots, n),$$
respectively, with a common namely positive factor k so that
$$H_i = kH_i' \ (i = 1, 2, 3, \ldots, n).$$
4. Explicitly determine the value of parameter k as the only unknown via this necessary and sufficient integral normalization condition
$$\Sigma_{i=1}^n H_i(c_{i+1} - c_{i-1}) = 2$$
so that
$$k = 2 / \Sigma_{i=0}^n H_i'(c_{i+1} - c_{i-1}).$$
5. Explicitly determine
$$H_i = kH_i' \ (i = 1, 2, 3, \ldots, n).$$

## 3.4. Mean Value (Mathematical Expectation)

Take the common integral definition [Cramér] of the mean value (mathematical expectation)
$$\mu = E(X) = \int_{-\infty}^{+\infty} xf(x)dx .$$
Use the corresponding formula for a piecewise linear probability distribution. Then in our continuous case we determine
$$\mu = \int_{-\infty}^{+\infty} xf(x)dx = \int_a^b xf(x)dx$$
$$= \Sigma_{i=0}^n (c_{i+1} - c_i)[R_i(2c_i + c_{i+1}) + L_{i+1}(c_i + 2c_{i+1})]/6$$
$$= 1/6 \; \Sigma_{i=0}^n (c_{i+1} - c_i)[H_i(c_{i+1} + 2c_i) + H_{i+1}(2c_{i+1} + c_i)]$$
$$= 1/6 \; \Sigma_{i=0}^n (c_{i+1} - c_i)H_i(c_{i+1} + 2c_i) + 1/6 \; \Sigma_{i=0}^n (c_{i+1} - c_i)H_{i+1}(2c_{i+1} + c_i)$$
$$= 1/6 \; \Sigma_{i=1}^n (c_{i+1} - c_i)H_i(c_{i+1} + 2c_i) + 1/6 \; \Sigma_{i=1}^n (c_i - c_{i-1})H_i(2c_i + c_{i-1})$$
$$= 1/6 \; \Sigma_{i=1}^n H_i[(c_{i+1} - c_i)(c_{i+1} + 2c_i) + (c_i - c_{i-1})(2c_i + c_{i-1})]$$
$$= 1/6 \; \Sigma_{i=1}^n H_i(c_{i+1}^2 + c_{i+1}c_i - 2c_i^2 + 2c_i^2 - c_ic_{i-1} - c_{i-1}^2)$$
$$= 1/6 \; \Sigma_{i=1}^n H_i(c_{i+1}^2 + c_{i+1}c_i - c_ic_{i-1} - c_{i-1}^2)$$
$$= 1/6 \; \Sigma_{i=1}^n H_i(c_{i+1} - c_{i-1})(c_{i+1} + c_i + c_{i-1})$$
and, finally,
$$\mu = \Sigma_{i=1}^n H_i(c_{i+1} - c_{i-1})(c_{i+1} + c_i + c_{i-1})/6.$$

## 3.5. Median Values

Use the common integral definition [Cramér] of median values $v$ for any of which both
$$P(X \leq v) \geq 1/2$$
and
$$P(X \geq v) \geq 1/2.$$
For a continual real-number random variable $X$,
$$P(X \leq v) = \int_{-\infty}^{v} f(x)dx = P(X \geq v) = \int_{v}^{+\infty} f(x)dx = 1/2.$$
To determine the set of all the median values $v$, we can use the same natural idea, way, and algorithm as for a general one-dimensional piecewise linear probability distribution but, naturally, with the formulas for a general one-dimensional piecewise linear continuous probability distribution.

## 3.6. Mode Values

To begin with, consider the common definition [Cramér] of mode values for any of which probability density distribution function f(x) takes its maximum value $f_{max}$ .

In particular, for a piecewise linear continuous probability distribution with probability density function f(x),

$$f_{max} = \max\{f(c_i) \mid i = 1, 2, \ldots, n\}.$$

If $f(x) = f_{max}$ at some x , then this x is one of the modes.

In particular, if $f(c_i) = f_{max}$ at some i , then $c_i$ at this i is one of the modes.

Nota bene: The set of all the modes both contains separate points

$$c_i \ (i = 1, 2, \ldots, n)$$

for which

$$f(c_i) = f_{sup} = f_{max}$$

and includes closed intervals

$$c_i \leq x \leq c_{i+1} \ (i = 1, 2, \ldots, n - 1)$$

for which

$$f(c_i) = f(c_{i+1}) = f_{sup} = f_{max} \ .$$

## 3.7. Variance

Take the common integral definition [Cramér] of the variance $\sigma^2$ of a random variable X as its second central moment, namely the squared standard deviation $\sigma$, or the expected value of the squared deviation from the mean:

$$\sigma^2 = E[(X - \mu)^2] = \int_{-\infty}^{+\infty} (x - \mu)^2 f(x) dx = \int_a^b (x - \mu)^2 f(x) dx.$$

Use the corresponding formula for a piecewise linear probability distribution. Then in our continuous case we determine

$$\begin{aligned}
\sigma^2 &= \Sigma_{i=0}^n (c_{i+1} - c_i)[R_i(c_{i+1}^2 + 2c_{i+1}c_i + 3c_i^2 - 4\mu c_{i+1} - 8\mu c_i + 6\mu^2) \\
&\quad + L_{i+1}(3c_{i+1}^2 + 2c_{i+1}c_i + c_i^2 - 8\mu c_{i+1} - 4\mu c_i + 6\mu^2)]/12 \\
&= \Sigma_{i=0}^n (c_{i+1} - c_i)[H_i(c_{i+1}^2 + 2c_{i+1}c_i + 3c_i^2 - 4\mu c_{i+1} - 8\mu c_i + 6\mu^2) \\
&\quad + H_{i+1}(3c_{i+1}^2 + 2c_{i+1}c_i + c_i^2 - 8\mu c_{i+1} - 4\mu c_i + 6\mu^2)]/12 \\
&= \Sigma_{i=0}^n (c_{i+1} - c_i)H_i(c_{i+1}^2 + 2c_{i+1}c_i + 3c_i^2 - 4\mu c_{i+1} - 8\mu c_i + 6\mu^2)/12 \\
&\quad + \Sigma_{i=0}^n (c_{i+1} - c_i)H_{i+1}(3c_{i+1}^2 + 2c_{i+1}c_i + c_i^2 - 8\mu c_{i+1} - 4\mu c_i + 6\mu^2)/12 \\
&= \Sigma_{i=1}^n (c_{i+1} - c_i)H_i(c_{i+1}^2 + 2c_{i+1}c_i + 3c_i^2 - 4\mu c_{i+1} - 8\mu c_i + 6\mu^2)/12 \\
&\quad + \Sigma_{i=0}^{n-1} (c_{i+1} - c_i)H_{i+1}(3c_{i+1}^2 + 2c_{i+1}c_i + c_i^2 - 8\mu c_{i+1} - 4\mu c_i + 6\mu^2)/12 \\
&= \Sigma_{i=1}^n (c_{i+1} - c_i)H_i(c_{i+1}^2 + 2c_{i+1}c_i + 3c_i^2 - 4\mu c_{i+1} - 8\mu c_i + 6\mu^2)/12 \\
&\quad + \Sigma_{i=1}^n (c_i - c_{i-1})H_i(3c_i^2 + 2c_i c_{i-1} + c_{i-1}^2 - 8\mu c_i - 4\mu c_{i-1} + 6\mu^2)/12 \\
&= \Sigma_{i=1}^n H_i(c_{i+1}^3 + 2c_{i+1}^2 c_i + 3c_{i+1}c_i^2 - 4\mu c_{i+1}^2 - 8\mu c_{i+1}c_i + 6\mu^2 c_{i+1} \\
&\quad - c_{i+1}^2 c_i - 2c_{i+1}c_i^2 - 3c_i^3 + 4\mu c_{i+1}c_i + 8\mu c_i^2 - 6\mu^2 c_i \\
&\quad + 3c_i^3 + 2c_i^2 c_{i-1} + c_i c_{i-1}^2 - 8\mu c_i^2 - 4\mu c_i c_{i-1} + 6\mu^2 c_i \\
&\quad - 3c_i^2 c_{i-1} - 2c_i c_{i-1}^2 - c_{i-1}^3 + 8\mu c_i c_{i-1} + 4\mu c_{i-1}^2 - 6\mu^2 c_{i-1})/12 \\
&= \Sigma_{i=1}^n H_i(c_{i+1}^3 + c_{i+1}^2 c_i + c_{i+1}c_i^2 - c_i^2 c_{i-1} - c_i c_{i-1}^2 - c_{i-1}^3 \\
&\quad - 4\mu c_{i+1}^2 - 4\mu c_{i+1}c_i + 4\mu c_i c_{i-1} + 4\mu c_{i-1}^2 + 6\mu^2 c_{i+1} - 6\mu^2 c_{i-1})/12 \\
&= \Sigma_{i=1}^n H_i(c_{i+1} - c_{i-1})(c_{i+1}^2 + c_{i+1}c_{i-1} + c_{i-1}^2 + c_{i+1}c_i + c_i c_{i-1} + c_i^2 \\
&\quad - 4\mu c_{i+1} - 4\mu c_{i-1} - 4\mu c_i + 6\mu^2)/12
\end{aligned}$$

and, finally,

$$\sigma^2 = \Sigma_{i=1}^n H_i(c_{i+1} - c_{i-1})[c_{i+1}^2 + c_i^2 + c_{i-1}^2 + c_{i+1}c_i + c_{i+1}c_{i-1} + c_i c_{i-1} - 4\mu(c_{i+1} + c_i + c_{i-1}) + 6\mu^2]/12$$

where

$$\mu = \Sigma_{i=1}^n H_i(c_{i+1} - c_{i-1})(c_{i+1} + c_i + c_{i-1})/6.$$

Nota bene: Similarly, we can also determine further initial and central moments etc. [Cramér], e.g. skewness

$$\gamma_1 = E[(X - \mu)^3/\sigma^3]$$

and excess

$$\gamma_2 = E[(X - \mu)^4/\sigma^4] - 3.$$

# 4. Tetragonal Probability Distribution

## 4.1. Main Definitions

A tetragonal probability distribution (Fig. 3) is a particular case of a general one-dimensional piecewise linear continuous probability distribution for n = 2 and further of a general one-dimensional piecewise linear probability distribution. Therefore, directly apply the above formulas for a general one-dimensional piecewise linear continuous probability distribution to a tetragonal probability distribution.

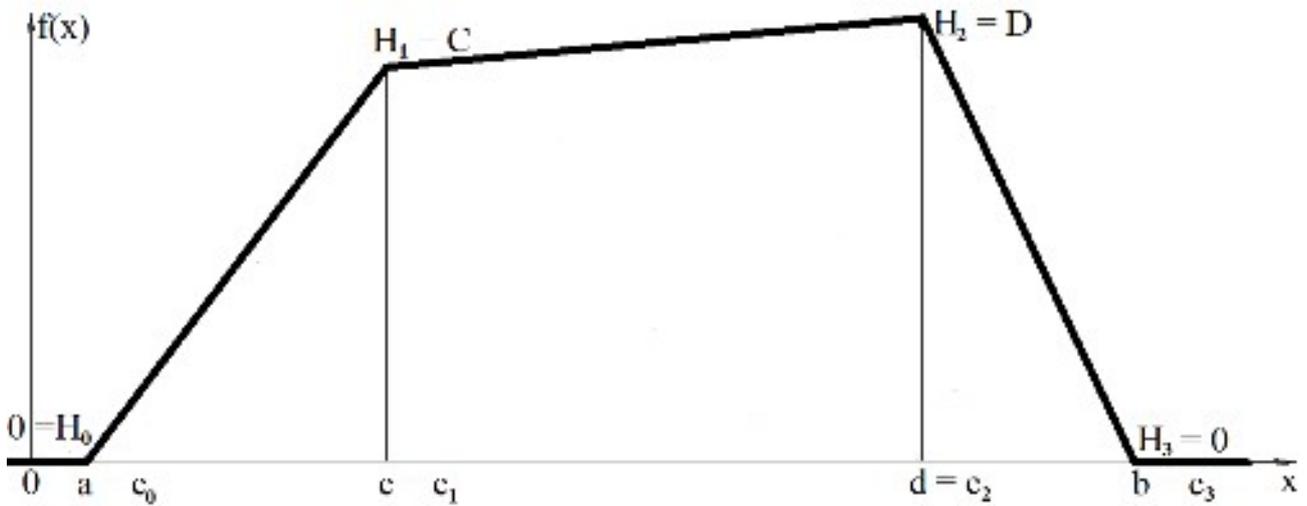

Figure 4. Tetragonal probability distribution

Here probability density distribution function f(x) is as always non-negative everywhere ($-\infty < x < +\infty$) and can be positive on some finite segment (closed interval)
$$-\infty < a \leq x \leq b < +\infty \ (a < b)$$
only. Let n = 2 intermediate points $c = c_1$ and $d = c_2$ in the non-decreasing order so that
$$a \leq c_1 \leq c_2 \leq b$$
divide this segment into n + 1 = 3 parts (pieces) of generally different lengths. To unify the notation, denote
$$c_0 = a,$$
$$c_3 = b,$$
$$c(i) = c_i \ (i = 0, 1, 2, 3).$$
On each of n + 1 = 3 closed intervals
$$c_i \leq x \leq c_{i+1} \ (i = 0, 1, 2),$$
probability density distribution function f(x) is linear. At n + 2 = 4 points
$$c_i \ (i = 0, 1, 2, 3),$$
f(x) takes finite non-negative values
$$H_i = f(c_i),$$
respectively. Naturally, we have
$$H_0 = 0,$$
$$H_3 = 0.$$
Note that

$$H_i = f(c_i) \ (i = 1, 2)$$

with additional natural notation

$$C = H_1,$$
$$D = H_2$$

for values f(x) at points

$$c = c_1,$$
$$d = c_2,$$

respectively, may be any finite non-negative values. At each of n + 2 = 4 points

$$c_i \ (i = 0, 1, 2, 3),$$

left and right one-sided limits

$$\lim f(x) = L_i \ (x \to c_i - 0),$$
$$\lim f(x) = R_i \ (x \to c_i + 0)$$

are equal to one another and coincide with $f(c_i)$. Therefore, we obtain

$$H_i = L_i = R_i \ (i = 0, 1, 2, 3).$$

Then on each of n + 1 = 3 closed intervals

$$c_i \leq x \leq c_{i+1} \ (i = 0, 1, 2),$$

linear probability density distribution function

$$f(x) = H_i + (H_{i+1} - H_i)(x - c_i)/(c_{i+1} - c_i)$$
$$= [H_i(c_{i+1} - x) + H_{i+1}(x - c_i)]/(c_{i+1} - c_i).$$

Use space-saving notation [Gelimson 2013] and the corresponding formula for a piecewise linear continuous probability distribution. Then in our case n = 2 we represent non-negative-valued function f(x) via

$$f_{[0,\infty)}(x_{(-\infty,\infty)}) = 0_{(-\infty, c(0)) \cup [c(3), \infty)} \cup \cup_{i=0}^{2} [H_i(c_{i+1} - x) + H_{i+1}(x - c_i)]/(c_{i+1} - c_i)_{[c(i), c(i+1))}$$

on the whole real axis $(-\infty, \infty)$.
Using

$$c_0 = a,$$
$$c_1 = c,$$
$$c_2 = d,$$
$$c_3 = b,$$
$$H_0 = 0,$$
$$H_1 = C,$$
$$H_2 = D,$$
$$H_3 = 0,$$

we obtain

$$f_{[0,\infty)}(x_{(-\infty,\infty)}) = 0_{(-\infty, c(0)) \cup [c(3), \infty)} \cup \{[H_0(c_1 - x) + H_1(x - c_0)]/(c_1 - c_0)\}_{[c(0), c(1))}$$
$$\cup \{[H_1(c_2 - x) + H_2(x - c_1)]/(c_2 - c_1)\}_{[c(1), c(2))} \cup \{[H_2(c_3 - x) + H_3(x - c_2)]/(c_3 - c_2)\}_{[c(2), c(3))},$$
$$f_{[0,\infty)}(x_{(-\infty,\infty)}) = 0_{(-\infty, a) \cup [b, \infty)} \cup C(x - a)/(c - a)_{[a, c)} \cup [C(d - x) + D(x - c)]/(d - c)_{[c, d)} \cup D(b - x)/(b - d)_{[d, b)}.$$

Integral (cumulative) probability distribution function

$$F(x) = P(X \leq x) = \int_{-\infty}^{x} f(t)dt$$

is probability $P(X \leq x)$ that real-number random variable X takes a real-number value not greater than x.

## 4.2. Normalization Condition

The probability of the event that X takes any finite real value is namely 1 because this event is certain. This gives integral normalization condition
$$\int_{-\infty}^{+\infty} f(x)dx = 1.$$
Use the corresponding formula for a piecewise linear continuous probability distribution. Then in our case n = 2 we determine
$$1 = \int_{-\infty}^{+\infty} f(x)dx = \int_a^b f(x)dx$$
$$= \Sigma_{i=1}^n H_i(c_{i+1} - c_{i-1})/2$$
$$= \Sigma_{i=1}^2 H_i(c_{i+1} - c_{i-1})/2.$$
We can also obtain this result at once rather geometrically than analytically, namely via adding the areas of the n + 1 = 3 rectangular trapezoids, among them 2 rectangular triangles at the endpoints a and b .

Therefore, to provide a possible (an admissible) probability density distribution function, necessary and sufficient integral normalization condition
$$\Sigma_{i=1}^2 H_i(c_{i+1} - c_{i-1}) = 2$$
has to be satisfied.

Using
$$c_0 = a ,$$
$$c_1 = c ,$$
$$c_2 = d ,$$
$$c_3 = b ,$$
$$H_1 = C ,$$
$$H_2 = D ,$$
we obtain
$$H_1(c_2 - c_0) + H_2(c_3 - c_1) = C(d - a) + D(b - c)$$
and, finally,
$$C(d - a) + D(b - c) = 2.$$

## 4.3. Normalization Algorithm

Nota bene: The obtained normalization condition is one condition only for
$$n + (n + 2) = 2n + 2 = 6$$
unknowns
$$H_i \ (i = 1, 2),$$
$$c_i \ (i = 0, 1, 2, 3).$$
Additionally,
$$H_i \geq 0 \ (i = 1, 2),$$
$$c_0 \leq c_1 \leq c_2 \leq c_3 \ .$$
Generally, it is not possible to simply take any admissible values of
$$2n + 2 - 1 = 2n + 1 = 5$$
unknowns and then to determine the value of the remaining unknown via this condition because it can happen that this value is inadmissible.

A natural idea, way, and algorithm to avoid this difficulty are as follows:

1. Fix
$$c_0 \leq c_1 \leq c_2 \leq c_3 \ .$$

2. Take any
$$H_i' \geq 0 \ (i = 1, 2)$$
so that there is at least one namely positive number among these n = 2 non-negative numbers.

3. Let
$$H_i \ (i = 1, 2)$$
be proportional to
$$H_i' \geq 0 \ (i = 1, 2),$$
respectively, with a common namely positive factor k so that
$$H_i = kH_i' \ (i = 1, 2).$$

4. Explicitly determine the value of parameter k as the only unknown via this necessary and sufficient integral normalization condition
$$\Sigma_{i=1}^2 H_i(c_{i+1} - c_{i-1}) = 2$$
so that
$$k = 2 / \Sigma_{i=0}^2 H_i'(c_{i+1} - c_{i-1}).$$

5. Explicitly determine
$$H_i = kH_i' \ (i = 1, 2).$$

Using
$$c_0 = a \ ,$$
$$c_1 = c \ ,$$
$$c_2 = d \ ,$$
$$c_3 = b \ ,$$
$$H_1 = C \ ,$$
$$H_2 = D$$
and naturally denoting
$$H_1' = C' \ ,$$
$$H_2' = D' \ ,$$
we obtain the same algorithm in the following form:

1. Fix
$$a \leq c \leq d \leq b \ .$$

2. Take any
$$C' \geq 0,$$
$$D' \geq 0$$

so that there is at least one namely positive number among these n = 2 non-negative numbers.
3. Let C and D be proportional to C' and D', respectively, with a common namely positive factor k so that
$$C = kC',$$
$$D = kD'.$$
4. Explicitly determine the value of parameter k as the only unknown via this necessary and sufficient integral normalization condition
$$C(d - a) + D(b - c) = 2$$
so that
$$k = 2/[C'(d - a) + D'(b - c)].$$
5. Explicitly determine
$$C = kC',$$
$$D = kD'.$$
We may also modify the same algorithm as follows:
1. Fix
$$a \leq c \leq d \leq b.$$
2. Take any
$$C' \geq 0,$$
$$D' \geq 0$$
so that there is at least one namely positive number among these n = 2 non-negative numbers.
3. Divide C' and D' by C' + D' to provide namely the unit sum of n = 2 non-negative numbers
$$w = C'/(C' + D'),$$
$$1 - w = D'/(C' + D').$$
4. Let C and D be proportional to w and 1 - w, respectively, with a common namely positive factor k so that
$$C = kw,$$
$$D = k(1 - w).$$
5. Explicitly determine the value of parameter k as the only unknown via this necessary and sufficient integral normalization condition
$$C(d - a) + D(b - c) = 2$$
so that
$$k = 2/[w(d - a) + (1 - w)(b - c)].$$
6. Explicitly determine
$$C = 2w/[w(d - a) + (1 - w)(b - c)],$$
$$D = 2(1 - w)/[w(d - a) + (1 - w)(b - c)].$$

## 4.4. Mean Value (Mathematical Expectation)

Take the common integral definition [Cramér] of the mean value (mathematical expectation)
$$\mu = E(X) = \int_{-\infty}^{+\infty} xf(x)dx .$$
Use the corresponding formula for a piecewise linear continuous probability distribution. Then in our case n = 2 we determine
$$\mu = \int_{-\infty}^{+\infty} xf(x)dx = \int_a^b xf(x)dx$$
$$= \sum_{i=1}^n H_i(c_{i+1} - c_{i-1})(c_{i+1} + c_i + c_{i-1})/6$$
$$= \sum_{i=1}^2 H_i(c_{i+1} - c_{i-1})(c_{i+1} + c_i + c_{i-1})/6.$$

Using
$$c_0 = a ,$$
$$c_1 = c ,$$
$$c_2 = d ,$$
$$c_3 = b ,$$
$$H_1 = C ,$$
$$H_2 = D ,$$
we obtain the same formula in the following form:
$$\mu = [H_1(c_2 - c_0)(c_2 + c_1 + c_0) + H_2(c_3 - c_1)(c_3 + c_2 + c_1)]/6,$$
$$\mu = [C(d - a)(d + c + a) + D(b - c)(b + d + c)]/6,$$
$$\mu = [C(d - a)(a + c + d) + D(b - c)(b + c + d)]/6.$$
Using
$$C = 2w/[w(d - a) + (1 - w)(b - c)],$$
$$D = 2(1 - w)/[w(d - a) + (1 - w)(b - c)],$$
finally determine
$$\mu = 1/3 \ [w(d - a)(a + c + d) + (1 - w)(b - c)(b + c + d)]/[w(d - a) + (1 - w)(b - c)].$$
To compare this result with the corresponding formula [Dorp Kotz] for 1 - w > 0, denote
$$\alpha = w/(1 - w)$$
and obtain
$$\mu = 1/3 \ [\alpha(d - a)(a + c + d) + (b - c)(b + c + d)]/[\alpha(d - a) + b - c].$$

## 4.5. Median Values

Use the common integral definition [Cramér] of median values $v$ for any of which both
$$P(X \leq v) \geq 1/2$$
and
$$P(X \geq v) \geq 1/2.$$
For a continual real-number random variable $X$,
$$P(X \leq v) = \int_{-\infty}^{v} f(x)dx = P(X \geq v) = \int_{v}^{+\infty} f(x)dx = 1/2.$$
To determine the set of all the median values $v$, we can use the same natural idea, way, and algorithm as for a general one-dimensional piecewise linear probability distribution but, naturally, with the formulas for a tetragonal probability distribution.

But using $n = 2$, make the same natural idea, way, and algorithm much more explicit:

1. First determine both
$$F(c) = \int_{-\infty}^{c} f(x)dx = \int_{a}^{c} f(x)dx = \int_{a}^{c} C(x - a)/(c - a) \, dx$$
$$= C/(c - a) \int_{a}^{c}(x - a)dx = C/(c - a) \left[(c^2 - a^2)/2 - a(c - a)\right]$$
$$= C[(c + a)/2 - a] = C(c - a)/2$$
and
$$F(d) = 1 - \int_{d}^{+\infty} f(x)dx = 1 - \int_{d}^{b} f(x)dx = 1 - \int_{d}^{b} D(b - x)/(b - d) \, dx$$
$$= 1 - D/(b - d) \int_{d}^{b}(b - x)dx = 1 - D/(b - d) \left[b(b - d) - (b^2 - d^2)/2\right]$$
$$= 1 - D[b - (b + d)/2] = 1 - D(b - d)/2.$$

2. If
$$F(c) > 1/2,$$
or, equivalently,
$$C(c - a) > 1,$$
then there is the only median value $v$ strictly between $a$ and $c$ so that
$$F(v) = 1/2,$$
$$F(v) = \int_{-\infty}^{v} f(x)dx = \int_{a}^{v} f(x)dx = \int_{a}^{v} C(x - a)/(c - a) \, dx$$
$$= C/(c - a) \int_{a}^{v}(x - a)dx = C/(c - a) \left[(v^2 - a^2)/2 - a(v - a)\right]$$
$$= C/(c - a) (v - a)^2/2 = 1/2,$$
$$(v - a)^2 = (c - a)/C,$$
$$v = a + [(c - a)/C]^{1/2}.$$

3. If
$$F(c) = 1/2,$$
or, equivalently,
$$C(c - a) = 1,$$
then there is the only median value
$$v = c.$$

4. If
$$F(d) < 1/2,$$
or, equivalently,
$$1 - D(b - d)/2 < 1/2,$$
$$D(b - d) > 1,$$
then there is the only median value $v$ strictly between $d$ and $b$ so that
$$F(v) = 1/2,$$
$$F(v) = 1 - \int_{v}^{+\infty} f(x)dx = 1 - \int_{v}^{b} f(x)dx = 1 - \int_{v}^{b} D(b - x)/(b - d) \, dx$$
$$= 1 - D/(b - d) \int_{v}^{b}(b - x)dx = 1 - D/(b - d) \left[b(b - v) - (b^2 - v^2)/2\right]$$
$$= 1 - D/(b - d) (b - v)^2/2 = 1/2,$$
$$D/(b - d) (b - v)^2 = 1,$$
$$(b - v)^2 = (b - d)/D,$$

$$v = b - [(b - d)/D]^{1/2} .$$

5. If
$$F(d) = 1/2,$$
or, equivalently,
$$1 - D(b - d)/2 = 1/2,$$
$$D(b - d) = 1,$$
then there is the only median value
$$v = d .$$

6. Finally, if
$$F(c) < 1/2 < F(d),$$
or, equivalently,
$$C(c - a) < 1$$
and
$$D(b - d) < 1,$$
then there is the only median value $v$ strictly between $c$ and $d$ ($c < v < d$) because incremental distribution function $F(c)$ strictly monotonically increases on this interval $(c , d)$ so that
$$F(v) = 1/2,$$
$$F(v) = \int_{-\infty}^{v} f(x)dx = \int_{a}^{v} f(x)dx = \int_{a}^{c} f(x)dx + \int_{c}^{v} f(x)dx$$
$$= F(c) + \int_{c}^{v} [C(d - x) + D(x - c)]/(d - c) \, dx$$
$$= C(c - a)/2 + \{C[d(v - c) - (v^2 - c^2)/2] + D[(v^2 - c^2)/2 - c(v - c)]\}/(d - c)$$
$$= C(c - a)/2 + [(Cd - Dc)(v - c) + (D - C)(v^2 - c^2)/2]/(d - c) = 1/2,$$
$$C(c - a)(d - c) + 2(Cd - Dc)(v - c) + (D - C)(v^2 - c^2) = d - c ,$$
$$(D - C)v^2 + 2(Cd - Dc)v + C(c - a)(d - c) - 2(Cd - Dc)c - (D - C)c^2 + c - d = 0.$$

6.1. If $D = C$ and, naturally, positive, then
$$2C(d - c)v + C(c - a)(d - c) - 2C(d - c)c + c - d = 0,$$
$$2Cv = 1 + C(a + c),$$
$$v = 1/(2C) + (a + c)/2.$$
Directly moving from left to right, we also obtain the same result
$$v = c + [1/2 - C(c - a)/2]/C$$
at once. We have
$$v - c = 1/(2C) + (a - c)/2 > 0$$
because
$$C(c - a) < 1.$$
Directly moving from right to left, we obtain
$$v = d - [1/2 - C(b - d)/2]/C = - 1/(2C) + d + (b - d)/2 = (b + d)/2 - 1/(2C)$$
at once. We have
$$d - v = d + 1/(2C) - (b + d)/2 > 0$$
because
$$C(b - d) < 1.$$
To prove the equivalence of these both formulas
$$v = 1/(2C) + (a + c)/2$$
and
$$v = (b + d)/2 - 1/(2C)$$
for $v$ , note that
$$1/(2C) + (a + c)/2 = (b + d)/2 - 1/(2C)$$
because the normalization condition
$$C(c - a)/2 + C(d - c) + C(b - d)/2 = 1$$
gives
$$(b - a + d - c)/2 = 1/C .$$
6.2. If $D \neq C$ , then there is the only median value $v$ strictly between $c$ and $d$ ($c < v < d$) because incremental distribution function F(c) strictly monotonically increases on this interval (c , d) so that
$$F(v) = 1/2.$$
Hence quadratic equation
$$(D - C)v^2 + 2(Cd - Dc)v + C(c - a)(d - c) - 2(Cd - Dc)c - (D - C)c^2 + c - d = 0$$
in v has exactly one solution on this interval (c , d).

## 4.6. Mode Values

To begin with, consider the common definition [Cramér] of mode values for any of which probability density distribution function f(x) takes its maximum value $f_{max}$ .
If C = D and, naturally, positive, then there are two modes c and d .
If C > D , then there is the only mode c .
If C < D , then there is the only mode d .

## 4.7. Variance

Take the common integral definition [Cramér] of the variance $\sigma^2$ of a random variable X as its second central moment, namely the squared standard deviation $\sigma$, or the expected value of the squared deviation from the mean:
$$\sigma^2 = E[(X - \mu)^2] = \int_{-\infty}^{+\infty} (x - \mu)^2 f(x) dx \, .$$
Use the corresponding formula for a piecewise linear continuous probability distribution. Then in our case n = 2 we determine
$$\sigma^2 = \int_{-\infty}^{+\infty} (x - \mu)^2 f(x) dx = \int_a^b (x - \mu)^2 f(x) dx$$
$$= \Sigma_{i=1}^n H_i(c_{i+1} - c_{i-1})[c_{i+1}^2 + c_i^2 + c_{i-1}^2 + c_{i+1}c_i + c_{i+1}c_{i-1} + c_i c_{i-1} - 4\mu(c_{i+1} + c_i + c_{i-1}) + 6\mu^2]/12$$
$$= \Sigma_{i=1}^2 H_i(c_{i+1} - c_{i-1})[c_{i+1}^2 + c_i^2 + c_{i-1}^2 + c_{i+1}c_i + c_{i+1}c_{i-1} + c_i c_{i-1} - 4\mu(c_{i+1} + c_i + c_{i-1}) + 6\mu^2]/12$$
where
$$\mu = \Sigma_{i=1}^2 H_i(c_{i+1} - c_{i-1})(c_{i+1} + c_i + c_{i-1})/6.$$
Using
$$c_0 = a \, ,$$
$$c_1 = c \, ,$$
$$c_2 = d \, ,$$
$$c_3 = b \, ,$$
$$H_1 = C \, ,$$
$$H_2 = D \, ,$$
we obtain the same formulas in the following forms:
$$\mu = [H_1(c_2 - c_0)(c_2 + c_1 + c_0) + H_2(c_3 - c_1)(c_3 + c_2 + c_1)]/6,$$
$$\mu = [C(d - a)(a + c + d) + D(b - c)(b + c + d)]/6,$$
as well as
$$\sigma^2 = H_1(c_2 - c_0)[c_2^2 + c_1^2 + c_0^2 + c_2 c_1 + c_2 c_0 + c_1 c_0 - 4\mu(c_2 + c_1 + c_0) + 6\mu^2]$$
$$+ H_2(c_3 - c_1)[c_3^2 + c_2^2 + c_1^2 + c_3 c_2 + c_3 c_1 + c_2 c_1 - 4\mu(c_3 + c_2 + c_1) + 6\mu^2]\}/12$$
$$= \{C(d - a)[d^2 + c^2 + a^2 + dc + da + ca - 4\mu(d + c + a) + 6\mu^2]$$
$$+ D(b - c)[b^2 + d^2 + c^2 + bd + bc + dc - 4\mu(b + d + c) + 6\mu^2]\}/12.$$
Finally,
$$\sigma^2 = \{C(d - a)[a^2 + c^2 + d^2 + ac + ad + cd - 4\mu(a + c + d) + 6\mu^2]$$
$$+ D(b - c)[b^2 + c^2 + d^2 + bc + bd + cd - 4\mu(b + c + d) + 6\mu^2]\}/12.$$
Nota bene: Similarly, we can also determine further initial and central moments etc. [Cramér], e.g. skewness
$$\gamma_1 = E[(X - \mu)^3/\sigma^3]$$
and excess
$$\gamma_2 = E[(X - \mu)^4/\sigma^4] - 3.$$

# 5. Piecewise Linear Probability Distribution Formulas Verification via a Triangular Probability Distribution

## 5.1. Main Definitions

Verify formulas for a general one-dimensional piecewise linear probability distribution using formulas [Cramér, Kotz Dorp, Wikipedia Triangular distribution] for a triangular probability distribution as a particular case of a general one-dimensional piecewise linear continuous probability distribution for n = 1 and further of a general one-dimensional piecewise linear probability distribution. Therefore, directly apply the above formulas for a general one-dimensional piecewise linear continuous probability distribution (or, alternatively, for a tetragonal probability distribution) to a triangular probability distribution (Fig. 4).

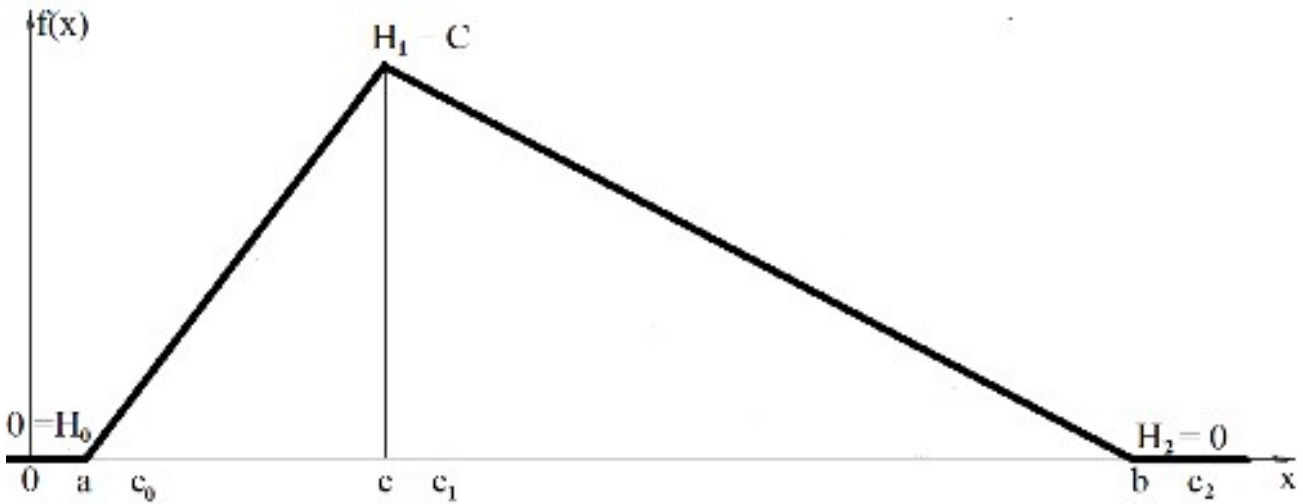

Figure 5. Triangular probability distribution

Here probability density distribution function f(x) is as always non-negative everywhere (-∞ < x < +∞) and can be positive on some finite segment (closed interval)
$$-\infty < a \leq x \leq b < +\infty \ (a < b)$$
only. Let n = 1 intermediate point $c = c_1$ so that
$$a \leq c_1 \leq b$$
divide this segment into n + 1 = 2 parts (pieces) of generally different lengths. To unify the notation, denote
$$c_0 = a ,$$
$$c_2 = b ,$$
$$c(i) = c_i \ (i = 0, 1, 2).$$
On each of n + 1 = 2 closed intervals
$$c_i \leq x \leq c_{i+1} \ (i = 0, 1),$$
probability density distribution function f(x) is linear. At n + 2 = 3 points
$$c_i \ (i = 0, 1, 2),$$
f(x) takes finite non-negative values
$$H_i = f(c_i),$$
respectively. Naturally, we have

$$H_0 = 0,$$
$$H_2 = 0.$$

Note that
$$H_1 = f(c_1)$$
with additional natural notation
$$C = H_1$$
for value f(x) at point
$$c = c_1$$
may be any finite positive value. At each of n + 2 = 3 points
$$c_i \ (i = 0, 1, 2),$$
left and right one-sided limits
$$\lim f(x) = L_i \ (x \to c_i - 0),$$
$$\lim f(x) = R_i \ (x \to c_i + 0)$$
are equal to one another and coincide with $f(c_i)$. Therefore, we obtain
$$H_i = L_i = R_i \ (i = 0, 1, 2).$$
Then on each of n + 1 = 2 closed intervals
$$c_i \leq x \leq c_{i+1} \ (i = 0, 1),$$
linear probability density distribution function
$$f(x) = H_i + (H_{i+1} - H_i)(x - c_i)/(c_{i+1} - c_i)$$
$$= [H_i(c_{i+1} - x) + H_{i+1}(x - c_i)]/(c_{i+1} - c_i).$$
Use space-saving notation [Gelimson 2013] and the corresponding formula for a piecewise linear continuous probability distribution. Then in our case n = 1 we represent non-negative-valued function f(x) via
$$f_{[0,\infty)}(x_{(-\infty,\infty)}) = 0_{(-\infty,c(0)) \cup [c(3),\infty)} \cup \cup_{i=0}^{1} [H_i(c_{i+1} - x) + H_{i+1}(x - c_i)]/(c_{i+1} - c_i)_{[c(i),c(i+1))}$$
on the whole real axis $(-\infty, \infty)$.
Using
$$c_0 = a,$$
$$c_1 = c,$$
$$c_2 = b,$$
$$H_0 = 0,$$
$$H_1 = C,$$
$$H_2 = 0,$$
we obtain
$$f_{[0,\infty)}(x_{(-\infty,\infty)}) = 0_{(-\infty,c(0)) \cup [c(2),\infty)} \cup [H_0(c_1 - x) + H_1(x - c_0)]/(c_1 - c_0)\}_{[c(0),c(1))}$$
$$\cup [H_1(c_2 - x) + H_2(x - c_1)]/(c_2 - c_1)\}_{[c(1),c(2))},$$
$$f_{[0,\infty)}(x_{(-\infty,\infty)}) = 0_{(-\infty,a) \cup [b,\infty)} \cup C(x - a)/(c - a)_{[a,c)} \cup C(b - x)/(b - c)_{[c,b)}.$$
Integral (cumulative) probability distribution function
$$F(x) = P(X \leq x) = \int_{-\infty}^{x} f(t)dt$$
is probability $P(X \leq x)$ that real-number random variable X takes a real-number value not greater than x .

## 5.2. Normalization Condition

The probability of the event that X takes any finite real value is namely 1 because this event is certain. This gives integral normalization condition
$$\int_{-\infty}^{+\infty} f(x)dx = 1.$$
Use the corresponding formula for a piecewise linear continuous probability distribution. Then in our case n = 1 we determine
$$1 = \int_{-\infty}^{+\infty} f(x)dx = \int_{a}^{b} f(x)dx$$
$$= \Sigma_{i=1}^{n} H_i(c_{i+1} - c_{i-1})/2 = \Sigma_{i=1}^{1} H_i(c_{i+1} - c_{i-1})/2 = H_1(c_2 - c_0).$$
We can also obtain this result at once rather geometrically than analytically, namely via adding the areas of the 2 rectangular triangles.

Therefore, to provide a possible (an admissible) probability density distribution function, necessary and sufficient integral normalization condition
$$H_1(c_2 - c_0) = 2$$
has to be satisfied.

Using
$$c_0 = a ,$$
$$c_1 = c ,$$
$$c_2 = b ,$$
$$H_1 = C ,$$
we obtain
$$H_1(c_2 - c_0) = C(b - a)$$
and, finally,
$$C(b - a) = 2,$$
$$C = 2/(b - a).$$
The known formulas [Cramér, Kotz Dorp, Wikipedia Triangular distribution] for a triangular probability distribution give the same result.

## 5.3. Mean Value (Mathematical Expectation)

Take the common integral definition [Cramér] of the mean value (mathematical expectation)
$$\mu = E(X) = \int_{-\infty}^{+\infty} xf(x)dx .$$
Use the corresponding formula for a piecewise linear continuous probability distribution. Then in our case n = 1 we determine
$$\mu = \int_{-\infty}^{+\infty} xf(x)dx = \int_a^b xf(x)dx$$
$$= \Sigma_{i=1}^n H_i(c_{i+1} - c_{i-1})(c_{i+1} + c_i + c_{i-1})/6$$
$$= \Sigma_{i=1}^1 H_i(c_{i+1} - c_{i-1})(c_{i+1} + c_i + c_{i-1})/6$$
and, finally,
$$\mu = H_1(c_2 - c_0)(c_2 + c_1 + c_0)/6.$$
Using
$$c_0 = a ,$$
$$c_1 = c ,$$
$$c_2 = b ,$$
$$H_1 = C ,$$
we obtain the same formula in the following form:
$$\mu = C(b - a)(b + c + a)/6.$$
Using
$$C = 2/(b - a),$$
finally obtain
$$\mu = (a + b + c)/3.$$
The known formulas [Kotz Dorp, Wikipedia Triangular distribution] for a triangular probability distribution give the same result.

# 5.4. Median Values

Use the common integral definition [Cramér] of median values $v$ for any of which both
$$P(X \leq v) \geq 1/2$$
and
$$P(X \geq v) \geq 1/2.$$
For a continual real-number random variable $X$,
$$P(X \leq v) = \int_{-\infty}^{v} f(x)dx = P(X \geq v) = \int_{v}^{+\infty} f(x)dx = 1/2.$$
To determine the set of all the median values $v$, we can use the same natural idea, way, and algorithm as for a general one-dimensional piecewise linear probability distribution but, naturally, with the formulas for a triangular probability distribution.
But using $n = 1$, as well as the corresponding algorithm and formulas for a tetragonal probability distribution with
$$d = c,$$
$$D = C,$$
$$C = 2/(b - a),$$
make the same natural idea, way, and algorithm as for a general one-dimensional piecewise linear probability distribution much more explicit:
1. First determine
$$F(c) = \int_{-\infty}^{c} f(x)dx = \int_{a}^{c} f(x)dx = \int_{a}^{c} C(x - a)/(c - a)\ dx$$
$$= C/(c - a) \int_{a}^{c} (x - a)dx = C/(c - a)\ [(c^2 - a^2)/2 - a(c - a)]$$
$$= C[(c + a)/2 - a] = C(c - a)/2 = (c - a)/(b - a).$$
2. If
$$F(c) > 1/2,$$
or, equivalently,
$$c > (a + b)/2,$$
then there is the only median value $v$ strictly between $a$ and $c$ so that
$$F(v) = 1/2,$$
$$F(v) = \int_{-\infty}^{v} f(x)dx = \int_{a}^{v} f(x)dx = \int_{a}^{v} C(x - a)/(c - a)\ dx$$
$$= C/(c - a) \int_{a}^{v} (x - a)dx = C/(c - a)\ [(v^2 - a^2)/2 - a(v - a)]$$
$$= C/(c - a)\ (v - a)^2/2 = 1/2,$$
$$(v - a)^2 = (c - a)/C,$$
$$v = a + [(c - a)/C]^{1/2},$$
$$v = a + [(b - a)(c - a)/2]^{1/2}.$$
The known formulas [Kotz Dorp, Wikipedia Triangular distribution] for a triangular probability distribution give the same result.
3. If
$$F(c) = 1/2,$$
or, equivalently,
$$c = (a + b)/2,$$
then there is the only median value
$$v = c = (a + b)/2.$$
Naturally, the known formulas [Cramér, Kotz Dorp, Wikipedia Triangular distribution] for a triangular probability distribution give the same obvious result.
4. If
$$F(c) < 1/2,$$
or, equivalently,
$$c < (a + b)/2,$$
then there is the only median value $v$ strictly between $c$ and $b$ so that

$$F(v) = 1/2,$$
$$F(v) = 1 - \int_v^{+\infty} f(x)dx = 1 - \int_v^b f(x)dx = 1 - \int_v^b C(b - x)/(b - c) \, dx$$
$$= 1 - C/(b - c) \int_v^b (b - x)dx = 1 - C/(b - c) [b(b - v) - (b^2 - v^2)/2]$$
$$= 1 - C/(b - c) (b - v)^2/2 = 1/2,$$
$$C/(b - c) (b - v)^2 = 1,$$
$$(b - v)^2 = (b - c)/C ,$$
$$v = b - [(b - c)/C]^{1/2}$$
$$v = b - [(b - a)(b - c)/2]^{1/2} .$$

The known formulas [Kotz Dorp, Wikipedia Triangular distribution] for a triangular probability distribution give the same result.

These three conditional formulas for the only median value $v$ can be unified as follows:
$$v = (a + b)/2 + \{[(b - a)(b - a + |2c - a - b|)]^{1/2} + a - b\}/2 \, \text{sign}(2c - a - b).$$

In fact, we obtain:

1) by $c > (a + b)/2$,
$$v = (a + b)/2 + \{[(b - a)(b - a + 2c - a - b)]^{1/2} + a - b\}/2$$
$$= (a + b)/2 + \{[(b - a)(2c - 2a)]^{1/2} + a - b\}/2$$
$$= a + [(b - a)(c - a)/2]^{1/2} ;$$

2) by $c = (a + b)/2$,
$$v = (a + b)/2;$$

3) by $c < (a + b)/2$,
$$v = (a + b)/2 - \{[(b - a)(b - a - 2c + a + b)]^{1/2} + a - b\}/2$$
$$= (a + b)/2 - \{[(b - a)(2b - 2c)]^{1/2} + a - b\}/2$$
$$= b - [(b - a)(b - c)/2]^{1/2} .$$

## 5.5. Mode Values

To begin with, consider the common definition [Cramér] of mode values for any of which probability density distribution function f(x) takes its maximum value $f_{max}$ .

In our case, there is the only mode c .

Naturally, the known formulas [Cramér, Kotz Dorp, Wikipedia Triangular distribution] for a triangular probability distribution give the same obvious result.

## 5.6. Variance

Take the common integral definition [Cramér] of the variance $\sigma^2$ of a random variable X as its second central moment, namely the squared standard deviation $\sigma$, or the expected value of the squared deviation from the mean:
$$\sigma^2 = E[(X - \mu)^2] = \int_{-\infty}^{+\infty} (x - \mu)^2 f(x)dx \ .$$
Use the corresponding formula for a piecewise linear continuous probability distribution. Then in our case n = 1 we determine
$$\sigma^2 = \int_{-\infty}^{+\infty} (x - \mu)^2 f(x)dx = \int_a^b (x - \mu)^2 f(x)dx$$
$$= \Sigma_{i=1}^n H_i(c_{i+1} - c_{i-1})[c_{i+1}^2 + c_i^2 + c_{i-1}^2 + c_{i+1}c_i + c_{i+1}c_{i-1} + c_ic_{i-1} - 4\mu(c_{i+1} + c_i + c_{i-1}) + 6\mu^2]/12$$
$$= \Sigma_{i=1}^1 H_i(c_{i+1} - c_{i-1})[c_{i+1}^2 + c_i^2 + c_{i-1}^2 + c_{i+1}c_i + c_{i+1}c_{i-1} + c_ic_{i-1} - 4\mu(c_{i+1} + c_i + c_{i-1}) + 6\mu^2]/12$$
and, finally,
$$\sigma^2 = H_1(c_2 - c_0)[c_2^2 + c_1^2 + c_0^2 + c_2c_1 + c_2c_0 + c_1c_0 - 4\mu(c_2 + c_1 + c_0) + 6\mu^2]/12$$
where
$$\mu = \Sigma_{i=1}^1 H_i(c_{i+1} - c_{i-1})(c_{i+1} + c_i + c_{i-1})/6 = H_1(c_2 - c_0)(c_2 + c_1 + c_0)/6.$$
Using
$$c_0 = a \ ,$$
$$c_1 = c \ ,$$
$$c_2 = b \ ,$$
$$H_1 = C \ ,$$
$$C = 2/(b - a),$$
or, alternatively, the above formulas for a tetragonal probability distribution with
$$d = c \ ,$$
$$D = C \ ,$$
we obtain the same formulas in the following forms:
$$\mu = C(b - a)(a + b + c)/6,$$
$$\mu = (a + b + c)/3,$$
as well as
$$\sigma^2 = C(b - a)[b^2 + c^2 + a^2 + bc + ba + ca - 4\mu(b + c + a) + 6\mu^2]/12,$$
$$\sigma^2 = [a^2 + b^2 + c^2 + ab + ac + bc - 4\mu(a + b + c) + 6\mu^2]/6,$$
Substituting
$$\mu = (a + b + c)/3,$$
we obtain
$$\sigma^2 = [a^2 + b^2 + c^2 + ab + ac + bc - 4/3 \ (a + b + c)^2 + 2/3 \ (a + b + c)^2]/6,$$
$$\sigma^2 = [3(a^2 + b^2 + c^2 + ab + ac + bc) - 2(a + b + c)^2]/18,$$
$$\sigma^2 = (3a^2 + 3b^2 + 3c^2 + 3ab + 3ac + 3bc - 2a^2 - 2b^2 - 2c^2 - 4ab - 4ac - 4bc)/18,$$
$$\sigma^2 = (a^2 + b^2 + c^2 - ab - ac - bc)/18.$$
Alternatively,
$$\sigma^2 = [(c - a)^2 + (b - c)^2 + (b - a)^2]/36.$$
The known formulas [Kotz & van Dorp, Wikipedia Triangular distribution] for a triangular probability distribution give the same result.
Nota bene: Similarly, we can also determine further initial and central moments etc. [Cramér], e.g. skewness
$$\gamma_1 = E[(X - \mu)^3/\sigma^3]$$
and excess
$$\gamma_2 = E[(X - \mu)^4/\sigma^4] - 3.$$

# Main Results and Conclusions

1. A piecewise linear probability distribution is very simple, natural, and typical, as well as sufficiently general.
2. A general one-dimensional piecewise linear probability distribution is very suitable for adequately modeling via efficiently approximating practically arbitrary nonlinear probability distribution with any required precision.
3. The explicit normalization, expectation, and variance formulas along with the median and mode formulas and algorithms for a general one-dimensional piecewise linear probability distribution are obtained and developed.
4. These formulas and algorithms are also applied to a general one-dimensional piecewise linear continuous probability distribution.
5. The formulas and algorithms for a general one-dimensional piecewise linear continuous probability distribution are very suitable for its important particular case, namely for a tetragonal probability distribution. It is also a natural generalization of a triangular probability distribution.
6. The known formulas for a triangular probability distribution as a further particular case of a general one-dimensional piecewise linear probability distribution provide verifying the obtained formulas and algorithms.
7. To additionally verify the present analytical methods, geometrical approach can be also applied if possible and useful.
8. The problems of the existence and uniqueness of the mean, median, and mode values for a general one-dimensional piecewise linear probability distribution are set and algorithmically solved.
9. The obtained formulas and developed algorithms have clear mathematical (probabilistic and statistical) sense and are simple and very suitable for setting and solving many typical urgent problems.
10. Piecewise linear probability distribution theory provides scientific basis for discovering and thoroughly investigating many complex phenomena and relations not only in probability theory and mathematical statistics, but also in physics, engineering, chemistry, biology, medicine, geology, astronomy, meteorology, agriculture, politics, management, economics, finance, psychology, etc.

# Bibliography

[Cramér] Harald Cramér. Mathematical Methods of Statistics. Princeton University Press, 1945

[van Dorp & Kotz] J. René van Dorp, Samuel Kotz. Generalized Trapezoidal Distributions. Metrika, Vol. 58, Issue 1, July 2003, 85-97. DOI:10.1007/s001840200230

[Encyclopaedia of Mathematics] Encyclopaedia of Mathematics. Ed. Michiel Hazewinkel. Volumes 1 to 10. Supplements I to III. Kluwer Academic Publ., Dordrecht, 1987-2002

[Gelimson 2003a] Lev Gelimson. Quantianalysis: Uninumbers, Quantioperations, Quantisets, and Multiquantities (now Uniquantities). Abhandlungen der WIGB (Wissenschaftlichen Gesellschaft zu Berlin), **3** (2003), Berlin, 15-21

[Gelimson 2003b] Lev Gelimson. General Problem Theory. Abhandlungen der WIGB (Wissenschaftlichen Gesellschaft zu Berlin), **3** (2003), Berlin, 26-32

[Gelimson 2012a] Lev Gelimson. Conditional Including Piecewise Functions, Functionals, Operators, Mappings, and Correspondences: Space-Saving Notation. Mathematical Journal of the "Collegium" All World Academy of Sciences, Munich (Germany), **12** (2012), 16

[Gelimson 2012b] Lev Gelimson. Inversion and Inverse: Notation, Discretization, and Continualization. Mathematical Journal of the "Collegium" All World Academy of Sciences, Munich (Germany), **12** (2012), 17

[Karlis & Xekalaki] Dimitris Karlis, Evdokia Xekalaki. The Polygonal Distribution. International Conference on Mathematical and Statistical Modeling in Honor of Enrique Castillo, June 28-30, 2006. Advances in Mathematical and Statistical Modeling Statistics for Industry and Technology 2008, pp 21-33. http://www.uclm.es/actividades0506/congresos/icmsm2006/articles/KarlisX06.pdf

[Kim] Jay J. Kim. Application of the Truncated Triangular and the Trapezoidal Distributions. Proceedings of the Survey Research Methods Section, American Statistical Association (2011). Session 191: Synthetic Data and Other Methods for Disclosure Limitation and Confidentiality Preservation. http://www.amstat.org/sections/srms/Proceedings/y2007/Files/JSM2007-000166.pdf

[Kotz & van Dorp] Samuel Kotz, Johan René van Dorp. Beyond Beta: Other Continuous Families of Distributions with Bounded Support and Applications. World Scientific, 2004. Chapter 1: The Triangular Distribution. http://books.google.de/books?id=NpK1xSB_SE4C http://www.worldscientific.com/doi/suppl/10.1142/5720/suppl_file/5720_chap1.pdf

[Wikipedia Triangular distribution] http://en.wikipedia.org/wiki/Triangular_distribution